\documentclass[submission]{FPSAC2022}

\usepackage[utf8x]{inputenc}
\usepackage[english]{babel}
\usepackage{amsmath,amsfonts,amssymb,amsthm}
\usepackage{mathtools}
\usepackage[T1]{fontenc}
\usepackage{hyperref}
\usepackage{cite}
\usepackage[inline]{enumitem}
\usepackage{dsfont}
\usepackage{multicol}
\usepackage{stmaryrd}
\usepackage{subfig}
\usepackage{multirow}
\usepackage{mleftright}
\usepackage{tikz}

\usepackage{xcolor}
\definecolor{ColA}{HTML}{7a644d}  % Marron
\definecolor{ColB}{HTML}{23be63} % Vert
\newcommand{\ColMix}[2]{\textcolor{ColA!#1!ColB}{#2}}

\newtheorem{Theorem}{Theorem}[section]
\newtheorem{Proposition}[Theorem]{Proposition}

\numberwithin{equation}{section}

\renewcommand{\geq}{\geqslant}

\title[Mockingbird lattices]{Mockingbird lattices}
\author[S. Giraudo]{%
    Samuele Giraudo%
    \thanks{
        \href{mailto:samuele.giraudo@univ-eiffel.fr}{\tt samuele.giraudo@univ-eiffel.fr}.
        This research has been partially supported by the projects CARPLO (ANR-20-CE40-0007)
        and LambdaComb (ANR-21-CE48-0017) of the Agence nationale de la recherche.}%
    \addressmark{1}}
\address{%
    \addressmark{1} LIGM, Univ. Gustave Eiffel, CNRS, ESIEE Paris, F-$77454$
    Marne-la-Vallée, France.}
\abstract{%
    We study combinatorial and order theoretic structures arising from the fragment of
    combinatory logic spanned by the basic combinator $\M$. This basic combinator, named as
    the Mockingbird by Smullyan, is defined by the rewrite rule $\M \VarX_1 \Rew \VarX_1
    \VarX_1$. We prove that the reflexive and transitive closure of this rewrite relation is
    a partial order on terms on $\M$ and that all connected components of its rewrite graph
    are Hasse diagrams of lattices. This last result is based on the introduction of
    lattices on some forests. We enumerate the elements, the edges of the Hasse diagrams,
    and the intervals of these lattices with the help of formal power series on terms and on
    forests.}
\keywords{%
    Partial orders; Lattices; Combinatory logic; Rewrite systems; Treelike structures;
    Formal power series.}
\received{\today}

\newcommand{\ColA}[1]{\textcolor{ColA}{#1}}

\newcommand{\Hide}[1]{\ColA{\tt HIDEN}}
\newcommand{\Def}[1]{\ColMix{50}{\em #1}}
\newcommand{\Par}[1]{\mleft(#1\mright)}
\newcommand{\Bra}[1]{\mleft\{#1\mright\}}
\newcommand{\Han}[1]{\mleft[#1\mright]}
\newcommand{\HanL}[1]{\mleft\llbracket#1\mright]}

\newcommand{\Angle}[1]{\mleft\langle#1\mright\rangle}
\newcommand{\AAngle}[1]{\Angle{\Angle{#1}}}

\newcommand{\OEIS}[1]{\href{http://oeis.org/#1}{{\bf #1}}}

\newcommand{\N}{\mathbb{N}}

\newcommand{\K}{\mathbb{K}}

\newcommand{\VarX}{\mathsf{x}}
\newcommand{\VarZ}{\mathsf{z}}

\newcommand{\TreeT}{\mathfrak{t}}
\newcommand{\TreeS}{\mathfrak{s}}
\newcommand{\TreeR}{\mathfrak{r}}

\newcommand{\ForestF}{\mathfrak{f}}
\newcommand{\ForestG}{\mathfrak{g}}

\newcommand{\SequenceA}{\mathbf{a}}
\newcommand{\SequenceB}{\mathbf{b}}

\newcommand{\SetVariables}{\mathbb{X}}
\newcommand{\GeneratingSet}{\mathfrak{G}}
\DeclareMathOperator{\Product}{\star}
\newcommand{\Combinator}[1]{\operatorname{\mathbf{#1}}}
\newcommand{\CLS}{\mathcal{C}}

\newcommand{\Rew}{\to}
\newcommand{\RewContext}{\Rightarrow}
\newcommand{\Leq}{\preccurlyeq}
\newcommand{\Equiv}{\equiv}
\newcommand{\SetTerms}{\mathfrak{T}}
\newcommand{\Deg}{\mathrm{deg}}
\newcommand{\Height}{\mathrm{ht}}
\newcommand{\RewGraph}{\mathrm{G}}
\newcommand{\SetDuplicativeTrees}{\mathcal{D}}
\newcommand{\SetDuplicativeForests}{\SetDuplicativeTrees^*}
\newcommand{\LeqDuplicative}{\ll}
\newcommand{\Meet}{\wedge}
\newcommand{\JJoin}{\vee}
\newcommand{\MockingbirdLattice}{\mathrm{M}}
\newcommand{\Poset}{\mathcal{P}}
\newcommand{\Ladder}{\mathfrak{l}}
\newcommand{\GreaterLadders}{\mathcal{L}}
\newcommand{\RightComb}{\mathfrak{r}}
\newcommand{\ToForest}{\texttt{fr}}
\newcommand{\T}{\mathrm{T}}
\newcommand{\M}{{\Combinator{M}}}
\newcommand{\C}{{\Combinator{C}}}
\newcommand{\Kk}{{\Combinator{K}}}
\newcommand{\Ss}{{\Combinator{S}}}
\newcommand{\Ii}{{\Combinator{I}}}
\newcommand{\Aa}{{\Combinator{A}}}
\newcommand{\Bb}{{\Combinator{B}}}
\newcommand{\Support}{\mathrm{Supp}}
\newcommand{\IndicatorFunction}{\mathds{1}}

\newcommand{\SeriesF}{\mathbf{f}}

\newcommand{\CharacteristicSeries}{\mathbf{c}}
\newcommand{\SeriesLadders}{\mathbf{ld}}

\newcommand{\EnumerationMap}{\mathrm{en}}
\newcommand{\SeriesGreater}{\mathbf{gr}}
\newcommand{\SeriesNbSmaller}{\mathbf{ns}}
\newcommand{\SeriesMeetDecomposition}{\mathbf{md}}
\newcommand{\SeriesCoverings}{\mathbf{cv}}
\newcommand{\SeriesNbInputs}{\mathbf{ni}}

\DeclareMathOperator{\HadamardProduct}{\boxtimes}
\DeclareMathOperator{\MaxProduct}{\uparrow}
\DeclareMathOperator{\ConcatenateForests}{\centerdot}
\newcommand{\Flat}{\mathrm{P}}
\newcommand{\Merge}{\mathrm{mg}}

\newcommand{\RewDuplicative}{%
    \!%
    \begin{tikzpicture}[Centering,xscale=.08]
        \node at(0,0){$\RewContext$};
        \node at(1,0){$\RewContext$};
    \end{tikzpicture}%
    \!}

\newcommand{\WhiteNode}{
    \begin{tikzpicture}[Centering]
        \node[Node,WhiteNode,minimum size=1.5mm]at(0,0){};
    \end{tikzpicture}}

\newcommand{\BlackNode}{
    \begin{tikzpicture}[Centering]
        \node[Node,BlackNode,minimum size=1.5mm]at(0,0){};
    \end{tikzpicture}}

\newcommand{\AnyNode}{
    \begin{tikzpicture}[Centering]
        \node[Node,AnyNode,minimum size=1.5mm]at(0,0){};
    \end{tikzpicture}}

\tikzstyle{Centering}=[{baseline={([yshift=-0.5ex]current bounding box.center)}}]
\tikzstyle{MarkA}=[draw=ColA!80,fill=ColA!8]
\tikzstyle{MarkB}=[draw=ColB!80,fill=ColB!8]

\tikzstyle{Node}=[circle,draw=ColB!90,fill=ColB!10,inner sep=1pt,minimum size=2mm,thick,
font=\scriptsize]
\tikzstyle{NodeST}=[font=\scriptsize]
\tikzstyle{LeafST}=[font=\scriptsize,node distance=2mm]
\tikzstyle{Edge}=[draw=ColB!95,cap=round,rounded corners=2.5pt,thick]
\tikzstyle{GraphVertex}=[circle,MarkA,inner sep=1pt,minimum size=1.5mm,thick]
\tikzstyle{GraphLabeledVertex}=[font=\footnotesize,node distance=3mm]
\tikzstyle{GraphEdge}=[ColB,cap=round,very thick]
\tikzstyle{GraphArc}=[GraphEdge,->]
\tikzstyle{WhiteNode}=[Node,draw=ColB!90,fill=ColB!4]
\tikzstyle{BlackNode}=[Node,draw=ColA!150,fill=ColA!70]
\tikzstyle{AnyNode}=[Node,rectangle,draw=ColA!150,fill=ColA!10]

\begin{document}
\maketitle

%%%%%%%%%%%%%%%%%%%%%%%%%%%%%%%%%%%%%%%%%%%%%%%%%%%%%%%%%%%%%%%%%%%%%%%%%%%%%%%%%%%%%%%%%%%%
%%%%%%%%%%%%%%%%%%%%%%%%%%%%%%%%%%%%%%%%%%%%%%%%%%%%%%%%%%%%%%%%%%%%%%%%%%%%%%%%%%%%%%%%%%%%
%%%%%%%%%%%%%%%%%%%%%%%%%%%%%%%%%%%%%%%%%%%%%%%%%%%%%%%%%%%%%%%%%%%%%%%%%%%%%%%%%%%%%%%%%%%%
\section*{Introduction}
Combinatory logic~\cite{HS08} is a model of computation introduced by
Schönfinkel~\cite{Sch24} and developed by Curry~\cite{Cur30} with the objective to abstain
from the need of bound variables specific to the $\lambda$-calculus. Its combinatorial heart
is formed by terms, which are binary trees with labeled leaves, and rules to compute a
result from a term, which are rewrite relations on trees~\cite{BKVT03}. An important
instance is the system containing the basic combinators $\Kk$ and $\Ss$ with the rewrite
rules $\Kk \VarX_1 \VarX_2 \Rew \VarX_1$ and $\Ss \VarX_1 \VarX_2 \VarX_3 \Rew \VarX_1
\VarX_3 \Par{\VarX_2 \VarX_3}$. This system is important because it is combinatorially
complete: each $\lambda$-term can be translated, by bracket abstraction
algorithms~\cite{Sch24} into a term over $\Kk$ and $\Ss$ emulating it.

A lot of other basic combinators with their own rewrite rules have been introduced by
Smullyan in~\cite{Smu85} after ---now widely used--- bird names, forming the enchanted
forest of combinator birds. For instance, $\Kk$ is the Kestrel and $\Ss$ is the Starling.
Usual computer science questions consist in considering a fragment of combinatory logic,
that is a finite set of basic combinators with their rewrite rules, and to ask whether
\begin{enumerate*}[label={\bf (\alph*)}]
    \item Given two terms $\TreeT$ and $\TreeT'$, can we decide if $\TreeT$ and $\TreeT'$
    can be rewritten eventually in a same term? This is known as the word
    problem~\cite{BKVT03,Sta00}. It admits a positive answer for some basic combinators like
    the Lark~\cite{Sta89,SWB93} and the Warbler~\cite{SWB93} but is still open for the
    Starling~\cite{BEJW17};
    \item Given a term $\TreeT$, can we decide if all rewrite sequences starting from
    $\TreeT$ are finite? This is known as the strong normalization problem. It admits a
    positive answer for the Starling~\cite{Wal00} and the Jay~\cite{PS01}.
\end{enumerate*}

Here, we pursue this study in a different direction by adopting a combinatorial, order
theoretic, and enumerative point of view. In particular, by denoting by $\Leq$ (resp.\
$\Equiv$) the reflexive and transitive (resp.\ reflexive, symmetric, and transitive) closure
of the rewrite relation, we try to
\begin{enumerate*}[label={\bf (\alph*')}]
    \item Determine if $\Leq$ is a partial order relation;
    \item Determine in this case if each interval of this poset is a lattice;
    \item Enumerate the $\Equiv$-equivalence classes of terms w.r.t.\ some size notions.
\end{enumerate*}
This work fits in this general project consisting in mixing combinatory logic with
combinatorics.

We start this project by studying the system made of the basic combinator $\M$, known as the
Mockingbird~\cite{Smu85,Sta17}. By drawing some portions of its rewrite graph, the first
properties that stand out are that the graph does not contain any nontrivial loops and that
its connected components are finite and have exactly one minimal and one maximal element. At
this stage, driven by computer exploration, we conjecture that $\Leq$ is a partial order
relation and that each $\Equiv$-equivalence class is a lattice w.r.t.\ $\Leq$. This lattice
property is for us a good clue for the fact that this system contains rather rich
combinatorial properties. To prove this, we introduce a new lattice on duplicative forests,
that are kinds of treelike structures, and construct a poset isomorphism between this last
poset and the poset on terms on $\M$. The Mockingbird lattice of order $d \geq 0$ is the
lattice $\MockingbirdLattice(d)$ consisting in the combinators on $\M$ greater than or equal
to the right comb combinator on $\M$ of degree $d$. Since any combinator on $\M$ can be seen
as a binary tree, this provides a new lattice structure on these objects. Many similar
lattices have been studied on binary trees such as among others the Tamari
lattice~\cite{Tam62}. However, unlike these lattices having for any $d \geq 0$ a cardinality
equal to the $d$-th Catalan number, the elements of $\MockingbirdLattice(d)$ are enumerated
by a different integer sequence. To obtain enumerative results about the Mockingbird
lattices and all the posets of terms on $\M$ in general, we use formal power series on terms
and on duplicative forests. In this way, we enumerate the maximal and minimal elements of
the poset of all terms on $\M$, and the cardinality, the number of edges of the Hasse
diagram, and the number of intervals of~$\MockingbirdLattice(d)$.

This paper is organized as follows. Section~\ref{sec:terms_rewrite_cls} contains definitions
about terms, rewrite relations, and combinatory logic systems. In
Section~\ref{sec:mockingbird_lattice}, we study the combinatory logic system on $\M$ and the
Mockingbird lattices. Section~\ref{sec:enumerative_properties} contains enumerative results.
This text ends with the presentation of some open questions.

%%%%%%%%%%%%%%%%%%%%%%%%%%%%%%%%%%%%%%%%%%%%%%%%%%%%%%%%%%%%%%%%%%%%%%%%%%%%%%%%%%%%%%%%%%%%
{\it General notations and conventions.}
For any integers $i$ and $j$, $[i, j]$ denotes the set $\{i, i + 1, \dots, j\}$. For any
integer $i$, $[i]$ denotes the set $[1, i]$ and $\HanL{i}$ denotes the set $[0, i]$. For any
set $A$, $A^*$ is the set of words on $A$. For any $w \in A^*$ and $a \in A$, $|w|_a$ is the
number of occurrences of $a$ in $w$. The only word of length $0$ is the empty word
$\epsilon$. If $P$ is a statement, we denote by $\IndicatorFunction_P$ the indicator
function (equals to $1$ if $P$ holds and $0$ otherwise).

%%%%%%%%%%%%%%%%%%%%%%%%%%%%%%%%%%%%%%%%%%%%%%%%%%%%%%%%%%%%%%%%%%%%%%%%%%%%%%%%%%%%%%%%%%%%
%%%%%%%%%%%%%%%%%%%%%%%%%%%%%%%%%%%%%%%%%%%%%%%%%%%%%%%%%%%%%%%%%%%%%%%%%%%%%%%%%%%%%%%%%%%%
%%%%%%%%%%%%%%%%%%%%%%%%%%%%%%%%%%%%%%%%%%%%%%%%%%%%%%%%%%%%%%%%%%%%%%%%%%%%%%%%%%%%%%%%%%%%
\section{Combinatory logic systems} \label{sec:terms_rewrite_cls}
An \Def{alphabet} is a finite set $\GeneratingSet$. Its elements are called \Def{basic
combinators}. Any element of the set $\SetVariables := \bigcup_{n \geq 1} \SetVariables_n$,
where $\SetVariables_n := \Bra{\VarX_1, \dots, \VarX_n}$, is a \Def{variable}. The set
$\SetTerms(\GeneratingSet)$ of \Def{$\GeneratingSet$-terms} is so that any variable of
$\SetVariables$ is a $\GeneratingSet$-term, any basic combinator of $\GeneratingSet$ is a
$\GeneratingSet$-term, and if $\TreeT_1$ and $\TreeT_2$ are two $\GeneratingSet$-terms, then
$\Par{\TreeT_1 \Product \TreeT_2}$ is a $\GeneratingSet$-term. Any term is thus a rooted
planar binary tree where leaves are decorated by variables or by basic combinators. We shall
express terms concisely by removing superfluous parentheses by considering that $\Product$
associates to the left and also by removing the symbols $\Product$. For instance, if
$\GeneratingSet = \Bra{\Aa, \Bb}$, the $\GeneratingSet$-term
\begin{math}
    \TreeT :=
    \Par{\Par{\Par{\Aa \Product \Bb} \Product \Par{\VarX_1 \Product \VarX_2}} \Product \Aa}
\end{math}
writes concisely as
\begin{math}
    \Aa \Bb \Par{\VarX_1 \VarX_2} \Aa.
\end{math}
Let $\TreeT$ be a $\GeneratingSet$-term. The \Def{degree} $\Deg(\TreeT)$ of $\TreeT$ is the
number of internal nodes of $\TreeT$ seen as a binary tree. The \Def{depth} of a node $u$ of
$\TreeT$ is the number of internal nodes in the path connecting the root of $\TreeT$ and
$u$. The \Def{height} $\Height(\TreeT)$ of $\TreeT$ is the maximal depth among all the nodes
of $\TreeT$. A \Def{combinator} is a term having no occurrence of any variable.

Let $\TreeT$ and $\TreeT'_1$, \dots, $\TreeT'_n$, $n \geq 0$, be $\GeneratingSet$-terms. The
\Def{composition} of $\TreeT$ with $\TreeT'_1$, \dots, $\TreeT'_n$ is the
$\GeneratingSet$-term $\TreeT \Han{\TreeT'_1, \dots, \TreeT'_n}$ obtained by simultaneously
replacing for all $i \in [n]$ all occurrences of the variables $\VarX_i$ in $\TreeT$
by~$\TreeT'_i$. For instance
\begin{math}
    \VarX_1 \Par{\Aa \VarX_1} \Par{\VarX_4 \VarX_2} \Han{\Bb, \VarX_1 \VarX_3}
    = \Bb \Par{\Aa \Bb} \Par{\VarX_4 \Par{\VarX_1 \VarX_3}}.
\end{math}
Given two $\GeneratingSet$-terms $\TreeT$ and $\TreeS$, $\TreeS$ is a \Def{factor} of
$\TreeT$ if
\begin{math}
    \TreeT
    = \TreeT' \Han{\TreeS_1, \dots, \TreeS_{i - 1}, \TreeS, \TreeS_{i + 1}, \dots, \TreeS_n}
    \Han{\TreeR_1, \dots, \TreeR_m}
\end{math}
for some integers $n, m \geq 0$ and $\GeneratingSet$-terms $\TreeT'$, $\TreeS_1$, \dots,
$\TreeS_{i - 1}$, $\TreeS_{i + 1}$, \dots, $\TreeS_n$, $\TreeR_1$, \dots, $\TreeR_m$,
where $\VarX_i$ appears in $\TreeT'$. When this property does not hold, $\TreeT$
\Def{avoids} $\TreeS$.

A \Def{rewrite relation} on $\SetTerms(\GeneratingSet)$ is a binary relation $\Rew$ on
$\SetTerms(\GeneratingSet)$. A \Def{combinatory logic system} (or \Def{CLS} for short) is a
pair $\CLS := (\GeneratingSet, \Rew)$ where $\GeneratingSet$ is an alphabet and $\Rew$ is a
rewrite relation on $\SetTerms(\GeneratingSet)$ such that for each basic combinator $\C$ of
$\GeneratingSet$, there is exactly one rule of the form
\begin{math}
    \C \VarX_1 \dots \VarX_n \Rew \TreeT_{\C}
\end{math}
where $n \geq 1$ and $\TreeT_{\C}$ is a term having no basic combinators and having all
variables in $\SetVariables_n$. The integer $n$ is the \Def{order} of $\C$ in $\CLS$. Some
well-known combinators $\C$ together with the terms $\TreeT_\C$ appearing among other
in~\cite{Smu85} are the Identity bird $\Ii$ of order $1$ with $\TreeT_{\Ii} = \VarX_1$, the
Mockingbird $\Combinator{M}$ of order $1$ with $\TreeT_\M = \VarX_1 \VarX_1$, the Kestrel
$\Kk$ of order $2$ with $\TreeT_{\Kk} = \VarX_1$, and the Starling $\Ss$ of order $3$ with
$\TreeT_{\Ss} = \VarX_1 \VarX_3 \Par{\VarX_2 \VarX_3}$. The \Def{context closure} of $\Rew$
is the binary relation $\RewContext$ on $\SetTerms(\GeneratingSet)$ defined as follows. For
any $\C \in \GeneratingSet$, by denoting by $n$ the order of $\C$, we have
\begin{math}
    \C \VarX_1 \dots \VarX_n \Han{\TreeS_1, \dots, \TreeS_n}
    \RewContext
    \TreeT_{\C} \Han{\TreeS_1, \dots, \TreeS_n}
\end{math}
for any $\TreeS_1, \dots, \TreeS_n \in \SetTerms(\GeneratingSet)$, and
\begin{math}
    \TreeT_1 \TreeT_2 \RewContext \TreeT'_1 \TreeT_2
\end{math}
for any $\TreeT_1, \TreeT_2 \in \SetTerms(\GeneratingSet)$ whenever $\TreeT_1 \RewContext
\TreeT'_1$, and
\begin{math}
    \TreeT_1 \TreeT_2 \RewContext \TreeT_1 \TreeT'_2
\end{math}
for any $\TreeT_1, \TreeT_2 \in \SetTerms(\GeneratingSet)$ whenever $\TreeT_2 \RewContext
\TreeT'_2$. For instance, if $\CLS$ is the CLS containing the basic combinators $\Kk$ and
$\Ss$, we have
\begin{math}
    \Ss (\Kk \Kk \Ss) \Kk (\Ss \Ss)
    \RewContext \Ss \Kk \Kk (\Ss \Ss)
    \RewContext \Kk (\Ss \Ss) (\Kk (\Ss \Ss))
    \RewContext \Ss \Ss.
\end{math}

Given a CLS $\CLS := (\GeneratingSet, \Rew)$, we denote by $\Leq$ the preorder defined as
the reflexive and transitive closure of $\RewContext$. The \Def{rewrite graph}
$\RewGraph_\CLS$ of $\CLS$ is the digraph of the binary relation $\RewContext$ on
$\SetTerms(\GeneratingSet)$. For any $\TreeT \in \SetTerms(\GeneratingSet)$,
$\RewGraph_\CLS(\TreeT)$ is the subgraph of $\RewGraph_\CLS$ restrained on
\begin{math}
    \Bra{\TreeT' \in \SetTerms(\GeneratingSet) : \TreeT \Leq \TreeT'}.
\end{math}
When $\Leq$ is antisymmetric, $\CLS$ has the \Def{poset property} and we denote by
$\Poset_\CLS$ the poset $(\SetTerms(\GeneratingSet), \Leq)$. For any $\TreeT \in
\SetTerms(\GeneratingSet)$, $\Poset_\CLS(\TreeT)$ is the subposet of $\Poset_\CLS$ having
$\TreeT$ as smallest element. When $\CLS$ has the poset property and, for any $\TreeT \in
\SetTerms(\GeneratingSet)$, $\Poset_\CLS(\TreeT)$ is a lattice, $\CLS$ has the \Def{lattice
property}. We denote by $\Equiv$ the equivalence relation defined as the reflexive,
symmetric, and transitive closure of $\RewContext$. If for any $\TreeT \in
\SetTerms(\GeneratingSet)$, the $\Equiv$-equivalence class $\Han{\TreeT}_\Equiv$ of $\TreeT$
is finite, then $\CLS$ is \Def{locally finite}. When $\CLS$ has the poset property and, for
any $\TreeT \in \SetTerms(\GeneratingSet)$, $\Han{\TreeT}_\Equiv$ has a unique minimal
element, $\CLS$ is \Def{rooted}. If for any $\TreeT, \TreeS_1, \TreeS_2 \in
\SetTerms(\GeneratingSet)$, $\TreeT \Leq \TreeS_1$ and $\TreeT \Leq \TreeS_2$ implies the
existence of $\TreeT' \in \SetTerms(\GeneratingSet)$ such that $\TreeS_1 \Leq  \TreeT'$ and
$\TreeS_2 \Leq \TreeT'$, then $\CLS$ is \Def{confluent}.

\begin{figure}[ht]
    \centering
    \subfloat[][$\RewGraph_\Rew(\Ii \Ii (\Ii \Ii \Ii))$.]{
    \begin{minipage}[c]{.14\textwidth}
        \centering
        \scalebox{.5}{
        \begin{tikzpicture}[Centering,xscale=1.2,yscale=1.3]
            \node[GraphLabeledVertex](0)at(0,0){$\Ii \Ii (\Ii \Ii \Ii)$};
            \node[GraphLabeledVertex](1)at(-1,-1){$\Ii (\Ii \Ii \Ii)$};
            \node[GraphLabeledVertex](2)at(1,-1){$\Ii \Ii (\Ii \Ii)$};
            \node[GraphLabeledVertex](3)at(-1,-2){$\Ii (\Ii \Ii)$};
            \node[GraphLabeledVertex](4)at(1,-2){$\Ii \Ii \Ii$};
            \node[GraphLabeledVertex](5)at(0,-3){$\Ii \Ii$};
            \node[GraphLabeledVertex](6)at(0,-4){$\Ii$};
            \draw[GraphArc](0)--(1);
            \draw[GraphArc](0)--(2);
            \draw[GraphArc](1)--(3);
            \draw[GraphArc](1)--(4);
            \draw[GraphArc](2)--(3);
            \draw[GraphArc](2)--(4);
            \draw[GraphArc](3)--(5);
            \draw[GraphArc](4)--(5);
            \draw[GraphArc](5)--(6);
        \end{tikzpicture}}
    \end{minipage}
    \label{subfig:rewrite_graph_I}}
    \hfill
    \subfloat[][$\RewGraph_\Rew(\Ss\Ss\Kk(\Ss\Ss)\Kk)$.]{
    \begin{minipage}[c]{.32\textwidth}
        \centering
        \scalebox{.5}{
        \begin{tikzpicture}[Centering,xscale=1.65,yscale=1.4]
            \node[GraphLabeledVertex](0)at(-4,0){$\Ss\Ss\Kk(\Ss\Ss)\Kk$};
            \node[GraphLabeledVertex](1)at(-4,-1){$\Ss(\Ss\Ss)(\Kk(\Ss\Ss))\Kk$};
            \node[GraphLabeledVertex](2)at(-4,-2){$\Ss\Ss\Kk(\Kk(\Ss\Ss)\Kk)$};
            \node[GraphLabeledVertex](3)at(-1,-3)
                {$\Ss(\Kk(\Ss\Ss)\Kk)(\Kk(\Kk(\Ss\Ss)\Kk))$};
            \node[GraphLabeledVertex](4)at(-4,-4){$\Ss\Ss\Kk(\Ss\Ss)$};
            \node[GraphLabeledVertex](5)at(-2,-4){$\Ss(\Ss\Ss)(\Kk(\Kk(\Ss\Ss)\Kk))$};
            \node[GraphLabeledVertex](6)at(0,-4){$\Ss(\Kk(\Ss\Ss)\Kk)(\Kk(\Ss\Ss))$};
            \node[GraphLabeledVertex](7)at(-1,-5){$\Ss(\Ss\Ss)(\Kk(\Ss\Ss))$};
            \draw[GraphArc](0)--(1);
            \draw[GraphArc](1)--(2);
            \draw[GraphArc](2)--(3);
            \draw[GraphArc](2)--(4);
            \draw[GraphArc](3)--(5);
            \draw[GraphArc](3)--(6);
            \draw[GraphArc](5)--(7);
            \draw[GraphArc](6)--(7);
            \draw[GraphArc](4)--(7);
        \end{tikzpicture}}
    \end{minipage}
    \label{subfig:rewrite_graph_KS}}
    \hfill
    \subfloat[][$\RewGraph_\Rew(\M (\M (\M \M)))$.]{
    \begin{minipage}[c]{.33\textwidth}
        \centering
        \scalebox{0.5}{
        \begin{tikzpicture}[Centering,xscale=0.9,yscale=0.9]
            \node[GraphLabeledVertex](M[M[MM]])at(-2,0){$\M(\M(\M\M))$};
            \node[GraphLabeledVertex](M[MM][M[MM]])at(4,-2){$\M(\M\M)(\M(\M\M))$};
            \node[GraphLabeledVertex](M[MM][MM[MM]])at(2,-4){$\M(\M\M)(\M\M(\M\M))$};
            \node[GraphLabeledVertex](MM[MM][M[MM]])at(6,-4){$\M\M(\M\M)(\M(\M\M))$};
            \node[GraphLabeledVertex](MM[MM][MM[MM]])at(4,-6){$\M\M(\M\M)(\M\M(\M\M))$};
            \node[GraphLabeledVertex](M[MM[MM]])at(-2,-4){$\M(\M\M(\M\M))$};
            \draw[GraphArc](M[M[MM]])--(M[MM][M[MM]]);
            \draw[GraphArc](M[MM][M[MM]])--(M[MM][MM[MM]]);
            \draw[GraphArc](M[MM][M[MM]])--(MM[MM][M[MM]]);
            \draw[GraphArc](M[MM][MM[MM]])--(MM[MM][MM[MM]]);
            \draw[GraphArc](MM[MM][M[MM]])--(MM[MM][MM[MM]]);
            \draw[GraphArc](M[M[MM]])--(M[MM[MM]]);
            \draw[GraphArc](M[MM[MM]])--(MM[MM][MM[MM]]);
            \draw[GraphArc](M[M[MM]])edge[loop left,looseness=4](M[M[MM]]);
            \draw[GraphArc](M[MM][M[MM]])edge[loop right,looseness=4](M[MM][M[MM]]);
            \draw[GraphArc](M[MM][MM[MM]])edge[loop left,looseness=4](M[MM][MM[MM]]);
            \draw[GraphArc](MM[MM][M[MM]])edge[loop right,looseness=4](MM[MM][M[MM]]);
            \draw[GraphArc](MM[MM][MM[MM]])edge[loop right,looseness=4](MM[MM][MM[MM]]);
            \draw[GraphArc](M[MM[MM]])edge[loop left,looseness=4](M[MM[MM]]);
        \end{tikzpicture}}
    \end{minipage}
    \label{subfig:rewrite_graph_M}}
    \caption{Some subgraphs of rewrite graphs of some CLS.}
    \label{fig:rewrite_graphs}
\end{figure}
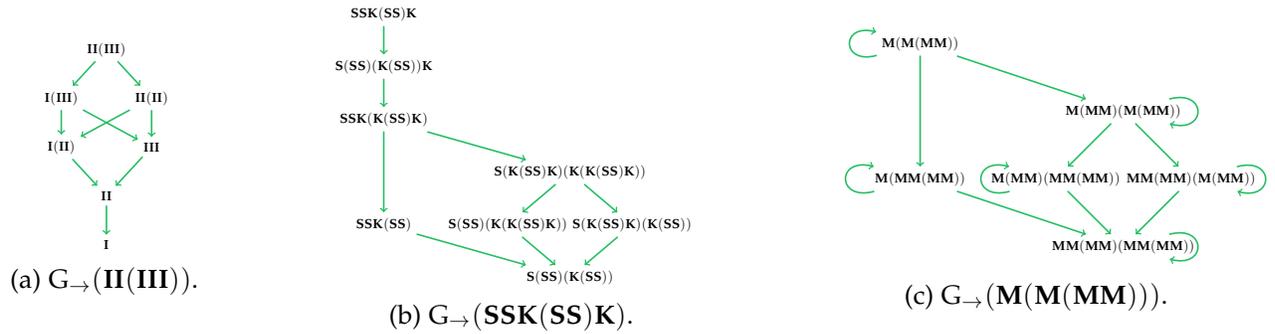

Consider for instance the CLS $\CLS$ containing the combinator $\Ii$. It is straightforward
to show that $\CLS$ has the poset property. Nevertheless, $\CLS$ has not the lattice
property, as suggested by the Hasse diagram shown in Figure~\ref{subfig:rewrite_graph_I}. It
is known that the CLS containing the combinators $\Kk$ and $\Ss$ has not the poset property.
Figure~\ref{subfig:rewrite_graph_KS} shows a subgraph of the rewrite graph of this CLS.
Figure~\ref{subfig:rewrite_graph_M} shows a subgraph of the rewrite graph of the CLS
containing the combinator $\M$. We shall study in details this CLS in the next sections.

A basic combinator $\C$ is \Def{hierarchical} if for any $i \in [n]$, $\VarX_i$ appears in
$\TreeT_\C$ at depth $n + 1 - i$. For instance, all the terms $\TreeT_\C$ such that $\C$ are
hierarchical and of order $3$ or less are $\VarX_1 \VarX_1$, $\VarX_1 \VarX_1 \VarX_2$,
$\VarX_2 \Par{\VarX_1 \VarX_1}$, $\VarX_1 \VarX_1 \VarX_2 \VarX_3$, $\VarX_2 \Par{\VarX_1
\VarX_1} \VarX_3$, $\VarX_3 \Par{\VarX_1 \VarX_1 \VarX_2}$, and $\VarX_3 \Par{\VarX_2
\Par{\VarX_1 \VarX_1}}$.

\begin{Proposition} \label{prop:confluence}
    Any CLS is confluent.
\end{Proposition}

\begin{Proposition} \label{prop:finite_equivalence_classes}
    If all basic combinators of a CLS $\CLS$ are hierarchical, then $\CLS$ is locally finite
    and all the $\GeneratingSet$-terms of a same connected component of $\RewGraph_\CLS$
    have the same height.
\end{Proposition}

Observe that by Propositions~\ref{prop:confluence}
and~\ref{prop:finite_equivalence_classes}, if $\CLS$ has the poset property and all its
basic combinators are hierarchical, then for any $\TreeT \in \SetTerms(\GeneratingSet)$, the
subposet $\Han{\TreeT}_\Equiv$ of $\Poset_\CLS$ has exactly one maximal element. If
additionally $\CLS$ is rooted, then for any $\TreeT \in \SetTerms(\GeneratingSet)$, the
subposet $\Han{\TreeT}_\Equiv$ of $\Poset_\CLS$ has exactly one minimal element.

%%%%%%%%%%%%%%%%%%%%%%%%%%%%%%%%%%%%%%%%%%%%%%%%%%%%%%%%%%%%%%%%%%%%%%%%%%%%%%%%%%%%%%%%%%%%
%%%%%%%%%%%%%%%%%%%%%%%%%%%%%%%%%%%%%%%%%%%%%%%%%%%%%%%%%%%%%%%%%%%%%%%%%%%%%%%%%%%%%%%%%%%%
%%%%%%%%%%%%%%%%%%%%%%%%%%%%%%%%%%%%%%%%%%%%%%%%%%%%%%%%%%%%%%%%%%%%%%%%%%%%%%%%%%%%%%%%%%%%
\section{The Mockingbird combinatory logic system} \label{sec:mockingbird_lattice}
Let $\CLS := (\GeneratingSet, \Rew)$ be the CLS such that $\GeneratingSet := \{\M\}$. We
call $\CLS$ the \Def{Mockingbird CLS}. From now, we shall simply write $\RewGraph$ instead
of $\RewGraph_\CLS$.

\begin{Proposition} \label{prop:first_graph_properties}
    The Mockingbird CLS is locally finite, has the poset property, and is rooted.
\end{Proposition}

Proposition~\ref{prop:first_graph_properties} is a consequence of the fact that $\M$ is
hierarchical and of Proposition~\ref{prop:finite_equivalence_classes}. By
Proposition~\ref{prop:first_graph_properties}, $\Poset_\CLS$ is a well-defined poset. From
now, we shall simply write $\Poset$ instead of $\Poset_\CLS$.

\begin{Proposition} \label{prop:maximal_minimal}
   A combinator $\TreeT$ is a maximal (resp.\ minimal) element of $\Poset$ if and only if
   $\TreeT$ avoids $\M \Par{\VarX_1 \VarX_2}$ (resp.\ $\Par{\VarX_1 \VarX_2} \Par{\VarX_1
   \VarX_2}$).
\end{Proposition}

A \Def{duplicative tree} is a planar rooted tree such that each node is either a \Def{black
node} $\BlackNode$ or a \Def{white node} $\WhiteNode$. A \Def{duplicative forest} is a word
$\ForestF$ of duplicative trees. We denote by $\SetDuplicativeTrees$ (resp.\
$\SetDuplicativeForests$) the set of such trees (resp.\ forests). The \Def{height}
$\Height(\ForestF)$ of $\ForestF$ is the number of internal nodes in a longest path
following edges connecting a node to one of its child. Each expression using some
occurrences of $\AnyNode$ denotes the two expressions obtained by replacing simultaneously
all $\AnyNode$ either by $\WhiteNode$ or by $\BlackNode$. The \Def{grafting product} is the
operation $\AnyNode$ on $\SetDuplicativeForests$ such that for any $\ForestF \in
\SetDuplicativeForests$, $\AnyNode(\ForestF)$ is the duplicative tree obtained by grafting
the roots of the duplicative trees of $\ForestF$ on a common root node~$\AnyNode$. The
\Def{concatenation product} is the binary operation $\ConcatenateForests$ on
$\SetDuplicativeForests$ such that for any $\ForestF_1, \ForestF_2 \in
\SetDuplicativeForests$, $\ForestF_1 \ConcatenateForests \ForestF_2$ is the duplicative
forest made of the trees of $\ForestF_1$ and then of the trees of $\ForestF_2$.

\begin{figure}[ht]
    \centering
    \scalebox{.75}{
    \begin{tikzpicture}[Centering,xscale=0.22,yscale=0.28]
        \tikzstyle{GraphEdge2}=[GraphEdge,thick]
        \node(WWW)at(-12,0){
            \begin{tikzpicture}[Centering,xscale=0.18,yscale=0.26]
                \node[WhiteNode](1)at(0.00,-2.67){};
                \node[WhiteNode](3)at(2.00,-1.33){};
                \node[WhiteNode](0)at(0.00,-1.33){};
                \draw[Edge](1)--(0);
            \end{tikzpicture}};
        \node(WBW)at(-12,-8){
            \begin{tikzpicture}[Centering,xscale=0.18,yscale=0.26]
                \node[BlackNode](1)at(0.00,-2.67){};
                \node[WhiteNode](3)at(2.00,-1.33){};
                \node[WhiteNode](0)at(0.00,-1.33){};
                \draw[Edge](1)--(0);
            \end{tikzpicture}};
        \node(BWWW)at(0,-4){
            \begin{tikzpicture}[Centering,xscale=0.18,yscale=0.26]
                \node[WhiteNode](0)at(0.00,-2.67){};
                \node[WhiteNode](2)at(2.00,-2.67){};
                \node[WhiteNode](4)at(4.00,-1.33){};
                \node[BlackNode](1)at(1.00,-1.33){};
                \draw[Edge](0)--(1);
                \draw[Edge](2)--(1);
            \end{tikzpicture}};
        \node(BBWW)at(-4,-8){
            \begin{tikzpicture}[Centering,xscale=0.18,yscale=0.26]
                \node[BlackNode](0)at(0.00,-2.67){};
                \node[WhiteNode](2)at(2.00,-2.67){};
                \node[WhiteNode](4)at(4.00,-1.33){};
                \node[BlackNode](1)at(1.00,-1.33){};
                \draw[Edge](0)--(1);
                \draw[Edge](2)--(1);
            \end{tikzpicture}};
        \node(BWBW)at(4,-8){
            \begin{tikzpicture}[Centering,xscale=0.18,yscale=0.26]
                \node[WhiteNode](0)at(0.00,-2.67){};
                \node[BlackNode](2)at(2.00,-2.67){};
                \node[WhiteNode](4)at(4.00,-1.33){};
                \node[BlackNode](1)at(1.00,-1.33){};
                \draw[Edge](0)--(1);
                \draw[Edge](2)--(1);
            \end{tikzpicture}};
        \node(BBBW)at(0,-12){
            \begin{tikzpicture}[Centering,xscale=0.18,yscale=0.26]
                \node[BlackNode](0)at(0.00,-2.67){};
                \node[BlackNode](2)at(2.00,-2.67){};
                \node[WhiteNode](4)at(4.00,-1.33){};
                \node[BlackNode](1)at(1.00,-1.33){};
                \draw[Edge](0)--(1);
                \draw[Edge](2)--(1);
            \end{tikzpicture}};
        \draw[GraphEdge2](WWW)--(WBW);
        \draw[GraphEdge2](WWW)--(BWWW);
        \draw[GraphEdge2](WBW)--(BBBW);
        \draw[GraphEdge2](BWWW)--(BBWW);
        \draw[GraphEdge2](BWWW)--(BWBW);
        \draw[GraphEdge2](BBWW)--(BBBW);
        \draw[GraphEdge2](BWBW)--(BBBW);
        \node(WWB)at(14,-8){
            \begin{tikzpicture}[Centering,xscale=0.18,yscale=0.26]
                \node[WhiteNode](1)at(0.00,-2.67){};
                \node[BlackNode](3)at(2.00,-1.33){};
                \node[WhiteNode](0)at(0.00,-1.33){};
                \draw[Edge](1)--(0);
            \end{tikzpicture}};
        \node(WBB)at(14,-16){
            \begin{tikzpicture}[Centering,xscale=0.18,yscale=0.26]
                \node[BlackNode](1)at(0.00,-2.67){};
                \node[BlackNode](3)at(2.00,-1.33){};
                \node[WhiteNode](0)at(0.00,-1.33){};
                \draw[Edge](1)--(0);
            \end{tikzpicture}};
        \node(BWWB)at(26,-12){
            \begin{tikzpicture}[Centering,xscale=0.18,yscale=0.26]
                \node[WhiteNode](0)at(0.00,-2.67){};
                \node[WhiteNode](2)at(2.00,-2.67){};
                \node[BlackNode](4)at(4.00,-1.33){};
                \node[BlackNode](1)at(1.00,-1.33){};
                \draw[Edge](0)--(1);
                \draw[Edge](2)--(1);
            \end{tikzpicture}};
        \node(BBWB)at(22,-16){
            \begin{tikzpicture}[Centering,xscale=0.18,yscale=0.26]
                \node[BlackNode](0)at(0.00,-2.67){};
                \node[WhiteNode](2)at(2.00,-2.67){};
                \node[BlackNode](4)at(4.00,-1.33){};
                \node[BlackNode](1)at(1.00,-1.33){};
                \draw[Edge](0)--(1);
                \draw[Edge](2)--(1);
            \end{tikzpicture}};
        \node(BWBB)at(30,-16){
            \begin{tikzpicture}[Centering,xscale=0.18,yscale=0.26]
                \node[WhiteNode](0)at(0.00,-2.67){};
                \node[BlackNode](2)at(2.00,-2.67){};
                \node[BlackNode](4)at(4.00,-1.33){};
                \node[BlackNode](1)at(1.00,-1.33){};
                \draw[Edge](0)--(1);
                \draw[Edge](2)--(1);
            \end{tikzpicture}};
        \node(BBBB)at(26,-20){
            \begin{tikzpicture}[Centering,xscale=0.18,yscale=0.26]
                \node[BlackNode](0)at(0.00,-2.67){};
                \node[BlackNode](2)at(2.00,-2.67){};
                \node[BlackNode](4)at(4.00,-1.33){};
                \node[BlackNode](1)at(1.00,-1.33){};
                \draw[Edge](0)--(1);
                \draw[Edge](2)--(1);
            \end{tikzpicture}};
        \draw[GraphEdge2](WWB)--(WBB);
        \draw[GraphEdge2](WWB)--(BWWB);
        \draw[GraphEdge2](WBB)--(BBBB);
        \draw[GraphEdge2](BWWB)--(BBWB);
        \draw[GraphEdge2](BWWB)--(BWBB);
        \draw[GraphEdge2](BBWB)--(BBBB);
        \draw[GraphEdge2](BWBB)--(BBBB);
        \draw[GraphEdge2](WWW)edge[bend left=16](WWB);
        \draw[GraphEdge2](WBW)edge[bend right=16](WBB);
        \draw[GraphEdge2](BWWW)edge[bend left=16](BWWB);
        \draw[GraphEdge2](BBWW)edge[bend right=8](BBWB);
        \draw[GraphEdge2](BWBW)edge[bend left=8](BWBB);
        \draw[GraphEdge2](BBBW)edge[bend right=16](BBBB);
    \end{tikzpicture}}
    \caption{The Hasse diagram of a maximal interval of the poset of duplicative forests.}
    \label{fig:example_poset_duplicative_forests}
\end{figure}
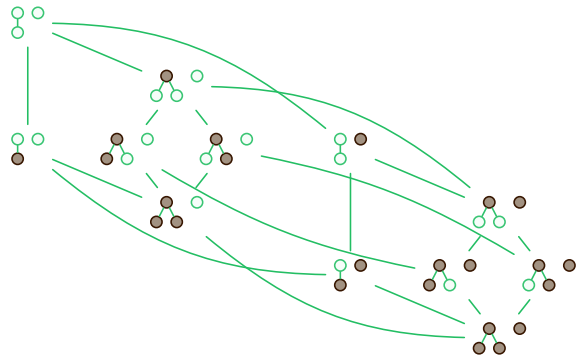

Let $\RewDuplicative$ be the binary relation on $\SetDuplicativeForests$ such that for any
$\ForestF, \ForestF' \in \SetDuplicativeForests$, we have $\ForestF \RewDuplicative
\ForestF'$ if $\ForestF'$ can be obtained from $\ForestF$ by selecting a white node of
$\ForestF$, by turning it into black, and by duplicating its sequence of descendants. For
instance, we have
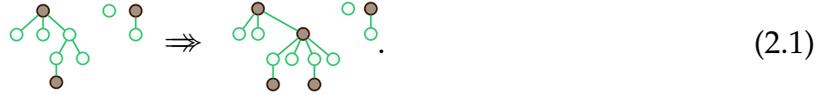
\begin{equation}
    \scalebox{.8}{
    \begin{tikzpicture}[Centering,xscale=0.22,yscale=0.18]
        \node[WhiteNode](0)at(-1.00,-4.40){};
        \node[WhiteNode](10)at(8.00,-4.40){};
        \node[WhiteNode](2)at(1.00,-4.40){};
        \node[BlackNode](4)at(2.00,-8.80){};
        \node[WhiteNode](6)at(4.00,-6.60){};
        \node[WhiteNode](8)at(6.00,-2.20){};
        \node[BlackNode](1)at(1.00,-2.20){};
        \node[WhiteNode](3)at(2.00,-6.60){};
        \node[WhiteNode](5)at(3.00,-4.40){};
        \node[BlackNode](9)at(8.00,-2.20){};
        \draw[Edge](0)--(1);
        \draw[Edge](10)--(9);
        \draw[Edge](2)--(1);
        \draw[Edge](3)--(5);
        \draw[Edge](4)--(3);
        \draw[Edge](5)--(1);
        \draw[Edge](6)--(5);
    \end{tikzpicture}}
    \enspace \RewDuplicative \enspace
    \scalebox{.8}{
    \begin{tikzpicture}[Centering,xscale=0.25,yscale=0.15]
        \node[WhiteNode](0)at(-0.25,-5.60){};
        \node[WhiteNode](11)at(7.00,-2.80){};
        \node[WhiteNode](13)at(8.50,-5.60){};
        \node[WhiteNode](2)at(1.00,-5.60){};
        \node[BlackNode](4)at(2.00,-11.20){};
        \node[WhiteNode](5)at(3.25,-8.40){};
        \node[BlackNode](8)at(4.75,-11.20){};
        \node[WhiteNode](9)at(6.00,-8.40){};
        \node[BlackNode](1)at(1.00,-2.80){};
        \node[BlackNode](12)at(8.50,-2.80){};
        \node[WhiteNode](3)at(2.00,-8.40){};
        \node[BlackNode](6)at(4.00,-5.60){};
        \node[WhiteNode](7)at(4.75,-8.40){};
        \draw[Edge](0)--(1);
        \draw[Edge](13)--(12);
        \draw[Edge](2)--(1);
        \draw[Edge](3)--(6);
        \draw[Edge](4)--(3);
        \draw[Edge](5)--(6);
        \draw[Edge](6)--(1);
        \draw[Edge](7)--(6);
        \draw[Edge](8)--(7);
        \draw[Edge](9)--(6);
    \end{tikzpicture}}.
\end{equation}
Observe that in this case, there are more black nodes in $\ForestF'$ than in $\ForestF$.
Hence, the reflexive and transitive closure $\LeqDuplicative$ of $\RewDuplicative$ is
antisymmetric so that $(\SetDuplicativeForests, \LeqDuplicative)$ is a poset. For any
$\ForestF \in \SetDuplicativeForests$, we denote by $\SetDuplicativeForests(\ForestF)$ the
subposet of $(\SetDuplicativeForests, \LeqDuplicative)$ on the set $\Bra{\ForestF' \in
\SetDuplicativeForests : \ForestF \LeqDuplicative \ForestF'}$.
Figure~\ref{fig:example_poset_duplicative_forests} shows the Hasse diagram of the poset
$\SetDuplicativeForests(\ForestF)$ for an $\ForestF \in \SetDuplicativeForests$.
According to this Hasse diagram, the poset $\SetDuplicativeForests$ is not graded.

Let $\Meet$ and $\JJoin$ be the two binary, commutative, and associative partial operations
on $\SetDuplicativeForests$ defined, for any $\ell \geq 0$, $\ForestF_1, \dots,
\ForestF_\ell \in \SetDuplicativeTrees$, $\ForestF_1', \dots, \ForestF_\ell' \in
\SetDuplicativeTrees$, and $\ForestF, \ForestF', \ForestF'' \in \SetDuplicativeForests$, by
\begin{minipage}{.47\textwidth}
\begin{equation} \begin{split}
    \ForestF_1 \dots \ForestF_\ell \Meet \ForestF_1' \dots \ForestF_\ell'
    & := \Par{\ForestF_1 \Meet \ForestF_1'} \dots \Par{\ForestF_\ell \Meet \ForestF_\ell'},
    \\
    \AnyNode(\ForestF) \Meet \AnyNode\Par{\ForestF'}
    & := \AnyNode\Par{\ForestF \Meet \ForestF'},
    \\
    \WhiteNode(\ForestF) \Meet \BlackNode\Par{\ForestF' \ForestF''}
    & := \WhiteNode\Par{\ForestF \Meet \ForestF' \Meet \ForestF''}
\end{split} \end{equation}
\end{minipage}
\hfill
\begin{minipage}{.47\textwidth}
\begin{equation} \begin{split}
    \ForestF_1 \dots \ForestF_\ell \JJoin \ForestF_1' \dots \ForestF_\ell'
    & := \Par{\ForestF_1 \JJoin \ForestF_1'}
        \dots \Par{\ForestF_\ell \JJoin \ForestF_\ell'},
    \\
    \AnyNode(\ForestF) \JJoin \AnyNode\Par{\ForestF'}
    & := \AnyNode\Par{\ForestF \JJoin \ForestF'},
    \\
    \WhiteNode(\ForestF) \JJoin \BlackNode\Par{\ForestF' \ForestF''}
    & := \BlackNode\Par{\Par{\ForestF \JJoin \ForestF'} \Par{\ForestF \JJoin \ForestF''}}.
\end{split} \end{equation}
\end{minipage}

\begin{Proposition} \label{prop:lattice_duplicative_forests}
    For any $\ForestF \in \SetDuplicativeForests$, the poset
    $\SetDuplicativeForests(\ForestF)$ is a lattice for the operations $\Meet$ and~$\JJoin$.
\end{Proposition}

We call \Def{duplicative forest lattice} of $\ForestF \in \SetDuplicativeForests$ the
lattice $\SetDuplicativeForests(\ForestF)$. To show that each subposet $\Poset(\TreeT)$,
$\TreeT \in \SetTerms(\GeneratingSet)$, of $\Poset$ is a lattice, we introduce a poset
isomorphism between $\Poset(\TreeT)$ and an interval of a lattice of duplicative forests.
For this, let $\ToForest : \SetTerms(\GeneratingSet) \to \SetDuplicativeForests$ be the map
defined, for any $\VarX_i \in \SetVariables$ and $\TreeT, \TreeT', \TreeT'' \in
\SetTerms(\GeneratingSet)$, by $\ToForest\Par{\VarX_i} := \ToForest(\M) := \ToForest\Par{\M
\M} := \epsilon$, $\ToForest\Par{\M \VarX_i} := \WhiteNode$, $\ToForest\Par{\M \Par{\TreeT
\TreeT'}} := \WhiteNode\Par{\ToForest\Par{\TreeT \TreeT'}}$, $\ToForest\Par{\VarX_i \TreeT}
:= \ToForest(\TreeT)$, and $\ToForest\Par{\Par{\TreeT \TreeT'} \TreeT''} :=
\ToForest\Par{\TreeT \TreeT'} \ConcatenateForests \ToForest\Par{\TreeT''}$. For instance,
\begin{equation}
    \scalebox{.8}{
    \begin{tikzpicture}[Centering,xscale=0.22,yscale=0.14]
        \node[LeafST](0)at(0.00,-11.57){$\M$};
        \node[LeafST](10)at(10.00,-23.14){$\VarX_1$};
        \node[LeafST](12)at(12.00,-11.57){$\VarX_1$};
        \node[LeafST](14)at(14.00,-15.43){$\M$};
        \node[LeafST](16)at(16.00,-15.43){$\M$};
        \node[LeafST](18)at(18.00,-11.57){$\M$};
        \node[LeafST](2)at(2.00,-19.29){$\M$};
        \node[LeafST](20)at(20.00,-15.43){$\M$};
        \node[LeafST](22)at(22.00,-15.43){$\VarX_3$};
        \node[LeafST](24)at(24.00,-11.57){$\VarX_2$};
        \node[LeafST](26)at(26.00,-11.57){$\VarX_2$};
        \node[LeafST](4)at(4.00,-19.29){$\VarX_2$};
        \node[LeafST](6)at(6.00,-19.29){$\M$};
        \node[LeafST](8)at(8.00,-23.14){$\M$};
        \node[NodeST](1)at(1.00,-7.71){$\Product$};
        \node[NodeST](11)at(11.00,-3.86){$\Product$};
        \node[NodeST](13)at(13.00,-7.71){$\Product$};
        \node[NodeST](15)at(15.00,-11.57){$\Product$};
        \node[NodeST](17)at(17.00,0.00){$\Product$};
        \node[NodeST](19)at(19.00,-7.71){$\Product$};
        \node[NodeST](21)at(21.00,-11.57){$\Product$};
        \node[NodeST](23)at(23.00,-3.86){$\Product$};
        \node[NodeST](25)at(25.00,-7.71){$\Product$};
        \node[NodeST](3)at(3.00,-15.43){$\Product$};
        \node[NodeST](5)at(5.00,-11.57){$\Product$};
        \node[NodeST](7)at(7.00,-15.43){$\Product$};
        \node[NodeST](9)at(9.00,-19.29){$\Product$};
        \draw[Edge](0)--(1);
        \draw[Edge](1)--(11);
        \draw[Edge](10)--(9);
        \draw[Edge](11)--(17);
        \draw[Edge](12)--(13);
        \draw[Edge](13)--(11);
        \draw[Edge](14)--(15);
        \draw[Edge](15)--(13);
        \draw[Edge](16)--(15);
        \draw[Edge](18)--(19);
        \draw[Edge](19)--(23);
        \draw[Edge](2)--(3);
        \draw[Edge](20)--(21);
        \draw[Edge](21)--(19);
        \draw[Edge](22)--(21);
        \draw[Edge](23)--(17);
        \draw[Edge](24)--(25);
        \draw[Edge](25)--(23);
        \draw[Edge](26)--(25);
        \draw[Edge](3)--(5);
        \draw[Edge](4)--(3);
        \draw[Edge](5)--(1);
        \draw[Edge](6)--(7);
        \draw[Edge](7)--(5);
        \draw[Edge](8)--(9);
        \draw[Edge](9)--(7);
        \node(r)at(17.00,3.25){};
        \draw[Edge](r)--(17);
    \end{tikzpicture}}
    \quad \xmapsto{\ToForest} \quad
    \scalebox{.9}{
    \begin{tikzpicture}[Centering,xscale=0.22,yscale=0.22]
        \node[WhiteNode](0)at(0.00,-3.50){};
        \node[WhiteNode](3)at(2.00,-5.25){};
        \node[WhiteNode](6)at(4.00,-3.50){};
        \node[WhiteNode](1)at(1.00,-1.75){};
        \node[WhiteNode](2)at(2.00,-3.50){};
        \node[WhiteNode](5)at(4.00,-1.75){};
        \draw[Edge](0)--(1);
        \draw[Edge](2)--(1);
        \draw[Edge](3)--(2);
        \draw[Edge](6)--(5);
    \end{tikzpicture}}.
\end{equation}
Immediately from the definition, we observe that this map is not injective. It can be shown
by structural induction on duplicative forests that the image of $\ToForest$ is the set of
all duplicative forests with no black nodes.

\begin{Proposition} \label{prop:poset_isomorphism}
    For any $\TreeT \in \SetTerms(\GeneratingSet)$, the posets $\Poset(\TreeT)$ and
    $\SetDuplicativeForests(\ToForest(\TreeT))$ are isomorphic.
\end{Proposition}

\begin{Theorem} \label{thm:mockingbird_lattices}
    For any $\TreeT \in \SetTerms(\GeneratingSet)$, the poset $\Poset(\TreeT)$ is a finite
    lattice.
\end{Theorem}

\begin{figure}[ht]
    \centering
    \subfloat[][$\MockingbirdLattice(0) \simeq \MockingbirdLattice(1)$.]{
    \begin{minipage}[b]{.16\textwidth}
        \centering
        \begin{tikzpicture}[Centering]
            \node[GraphVertex](1)at(0,0){};
        \end{tikzpicture}
    \end{minipage}
    \label{subfig:M_0_M_1}}
    \hfill
    \subfloat[][$\MockingbirdLattice(2)$.]{
    \begin{minipage}[b]{.11\textwidth}
        \centering
        \begin{tikzpicture}[Centering,yscale=.5]
            \tikzstyle{GraphEdge2}=[GraphEdge,thin]
            \node[GraphVertex](1)at(0,0){};
            \node[GraphVertex](2)at(0,-1){};
            \draw[GraphEdge2](1)--(2);
        \end{tikzpicture}
    \end{minipage}
    \label{subfig:M_2}}
    \hfill
    \subfloat[][$\MockingbirdLattice(3)$.]{
    \begin{minipage}[b]{.11\textwidth}
        \centering
        \begin{tikzpicture}[Centering,xscale=.22,yscale=.45]
            \tikzstyle{GraphEdge2}=[GraphEdge,thin]
            \node[GraphVertex](1)at(-1,0){};
            \node[GraphVertex](2)at(-1,-2){};
            \node[GraphVertex](3)at(2,-1){};
            \node[GraphVertex](4)at(1,-2){};
            \node[GraphVertex](5)at(3,-2){};
            \node[GraphVertex](6)at(2,-3){};
            \draw[GraphEdge2](1)--(2);
            \draw[GraphEdge2](1)--(3);
            \draw[GraphEdge2](3)--(4);
            \draw[GraphEdge2](3)--(5);
            \draw[GraphEdge2](4)--(6);
            \draw[GraphEdge2](5)--(6);
            \draw[GraphEdge2](2)--(6);
        \end{tikzpicture}
    \end{minipage}
    \label{subfig:M_3}}
    \hfill
    \subfloat[][$\MockingbirdLattice(4)$.]{
    \begin{minipage}[b]{.56\textwidth}
        \centering
        \begin{tikzpicture}[Centering,xscale=.08,yscale=.08]
            \tikzstyle{GraphVertex2}=[GraphVertex,minimum size=1mm]
            \tikzstyle{GraphEdge2}=[GraphEdge,thin]
            \node[GraphVertex2] (0) at (-28, -26) {};
            \node[GraphVertex2] (1) at (-22, -28) {};
            \node[GraphVertex2] (2) at (-24, -30) {};
            \node[GraphVertex2] (3) at (-20, -30) {};
            \node[GraphVertex2] (4) at (-22, -32) {};
            \node[GraphVertex2] (5) at (-28, -30) {};
            \node[GraphVertex2] (6) at (28, 0) {};
            \node[GraphVertex2] (7) at (34, -2) {};
            \node[GraphVertex2] (8) at (32, -4) {};
            \node[GraphVertex2] (9) at (36, -4) {};
            \node[GraphVertex2] (10) at (34, -6) {};
            \node[GraphVertex2] (11) at (28, -4) {};
            \node[GraphVertex2] (12) at (34, -20) {};
            \node[GraphVertex2] (13) at (40, -22) {};
            \node[GraphVertex2] (14) at (38, -24) {};
            \node[GraphVertex2] (15) at (42, -24) {};
            \node[GraphVertex2] (16) at (40, -26) {};
            \node[GraphVertex2] (17) at (34, -24) {};
            \node[GraphVertex2] (18) at (1, -14) {};
            \node[GraphVertex2] (19) at (7, -16) {};
            \node[GraphVertex2] (20) at (5, -18) {};
            \node[GraphVertex2] (21) at (9, -18) {};
            \node[GraphVertex2] (22) at (7, -20) {};
            \node[GraphVertex2] (23) at (1, -18) {};
            \node[GraphVertex2] (24) at (7, -34) {};
            \node[GraphVertex2] (25) at (13, -36) {};
            \node[GraphVertex2] (26) at (11, -38) {};
            \node[GraphVertex2] (27) at (15, -38) {};
            \node[GraphVertex2] (28) at (13, -40) {};
            \node[GraphVertex2] (29) at (7, -38) {};
            \node[GraphVertex2] (30) at (-7, 8) {};
            \node[GraphVertex2] (31) at (-1, 6) {};
            \node[GraphVertex2] (32) at (-3, 4) {};
            \node[GraphVertex2] (33) at (1, 4) {};
            \node[GraphVertex2] (34) at (-1, 2) {};
            \node[GraphVertex2] (35) at (-7, 4) {};
            \node[GraphVertex2] (36) at (-3, 18) {};
            \node[GraphVertex2] (37) at (38, 12) {};
            \node[GraphVertex2] (38) at (9, -4) {};
            \node[GraphVertex2] (39) at (46, -10) {};
            \node[GraphVertex2] (40) at (17, -26) {};
            \node[GraphVertex2] (41) at (-24, -16) {};
            \draw[GraphEdge2] (0)--(1);
            \draw[GraphEdge2] (1)--(2);
            \draw[GraphEdge2] (1)--(3);
            \draw[GraphEdge2] (3)--(4);
            \draw[GraphEdge2] (2)--(4);
            \draw[GraphEdge2] (0)--(5);
            \draw[GraphEdge2] (5)--(4);
            \draw[GraphEdge2] (6)--(7);
            \draw[GraphEdge2] (7)--(8);
            \draw[GraphEdge2] (7)--(9);
            \draw[GraphEdge2] (9)--(10);
            \draw[GraphEdge2] (8)--(10);
            \draw[GraphEdge2] (6)--(11);
            \draw[GraphEdge2] (11)--(10);
            \draw[GraphEdge2] (12)--(13);
            \draw[GraphEdge2] (13)--(14);
            \draw[GraphEdge2] (13)--(15);
            \draw[GraphEdge2] (15)--(16);
            \draw[GraphEdge2] (14)--(16);
            \draw[GraphEdge2] (12)--(17);
            \draw[GraphEdge2] (17)--(16);
            \draw[GraphEdge2] (18)--(19);
            \draw[GraphEdge2] (19)--(20);
            \draw[GraphEdge2] (19)--(21);
            \draw[GraphEdge2] (21)--(22);
            \draw[GraphEdge2] (20)--(22);
            \draw[GraphEdge2] (18)--(23);
            \draw[GraphEdge2] (23)--(22);
            \draw[GraphEdge2] (24)--(25);
            \draw[GraphEdge2] (25)--(26);
            \draw[GraphEdge2] (25)--(27);
            \draw[GraphEdge2] (27)--(28);
            \draw[GraphEdge2] (26)--(28);
            \draw[GraphEdge2] (24)--(29);
            \draw[GraphEdge2] (29)--(28);
            \draw[GraphEdge2] (30)--(31);
            \draw[GraphEdge2] (31)--(32);
            \draw[GraphEdge2] (31)--(33);
            \draw[GraphEdge2] (33)--(34);
            \draw[GraphEdge2] (32)--(34);
            \draw[GraphEdge2] (30)--(35);
            \draw[GraphEdge2] (35)--(34);
            \draw[GraphEdge2] (36)--(37);
            \draw[GraphEdge2] (37)--(38);
            \draw[GraphEdge2] (37)--(39);
            \draw[GraphEdge2] (39)--(40);
            \draw[GraphEdge2] (38)--(40);
            \draw[GraphEdge2] (36)--(41);
            \draw[GraphEdge2] (41)--(40);
            \draw[GraphEdge2] (30)--(0);
            \draw[GraphEdge2] (35)--(5);
            \draw[GraphEdge2] (31)--(1);
            \draw[GraphEdge2] (32)--(2);
            \draw[GraphEdge2] (33)--(3);
            \draw[GraphEdge2] (34)--(4);
            \draw[GraphEdge2] (30)--(6);
            \draw[GraphEdge2] (35)--(11);
            \draw[GraphEdge2] (32)--(8);
            \draw[GraphEdge2] (34)--(10);
            \draw[GraphEdge2] (6)--(18);
            \draw[GraphEdge2] (11)--(23);
            \draw[GraphEdge2] (8)--(20);
            \draw[GraphEdge2] (9)--(21);
            \draw[GraphEdge2] (7)--(19);
            \draw[GraphEdge2] (10)--(22);
            \draw[GraphEdge2] (6)--(12);
            \draw[GraphEdge2] (7)--(13);
            \draw[GraphEdge2] (11)--(17);
            \draw[GraphEdge2] (8)--(14);
            \draw[GraphEdge2] (9)--(15);
            \draw[GraphEdge2] (10)--(16);
            \draw[GraphEdge2] (0)--(24);
            \draw[GraphEdge2] (18)--(24);
            \draw[GraphEdge2] (12)--(24);
            \draw[GraphEdge2] (5)--(29);
            \draw[GraphEdge2] (23)--(29);
            \draw[GraphEdge2] (17)--(29);
            \draw[GraphEdge2] (2)--(26);
            \draw[GraphEdge2] (20)--(26);
            \draw[GraphEdge2] (14)--(26);
            \draw[GraphEdge2] (3)--(27);
            \draw[GraphEdge2] (21)--(27);
            \draw[GraphEdge2] (15)--(27);
            \draw[GraphEdge2] (1)--(25);
            \draw[GraphEdge2] (19)--(25);
            \draw[GraphEdge2] (13)--(25);
            \draw[GraphEdge2] (4)--(28);
            \draw[GraphEdge2] (22)--(28);
            \draw[GraphEdge2] (16)--(28);
            \draw[GraphEdge2] (36)--(30);
            \draw[GraphEdge2] (41)--(5);
            \draw[GraphEdge2] (37)--(7);
            \draw[GraphEdge2] (38)--(20);
            \draw[GraphEdge2] (39)--(15);
            \draw[GraphEdge2] (40)--(28);
            \draw[GraphEdge2] (31)--(7);
            \draw[GraphEdge2] (33)--(9);
        \end{tikzpicture}
    \end{minipage}
    \label{subfig:M_4}}
    \caption{The Hasse diagrams of the Mockingbird lattices $\MockingbirdLattice(d)$ for
    $d \in \HanL{4}$.}
    \label{fig:hasse_diagram_mockingbird_lattices}
\end{figure}
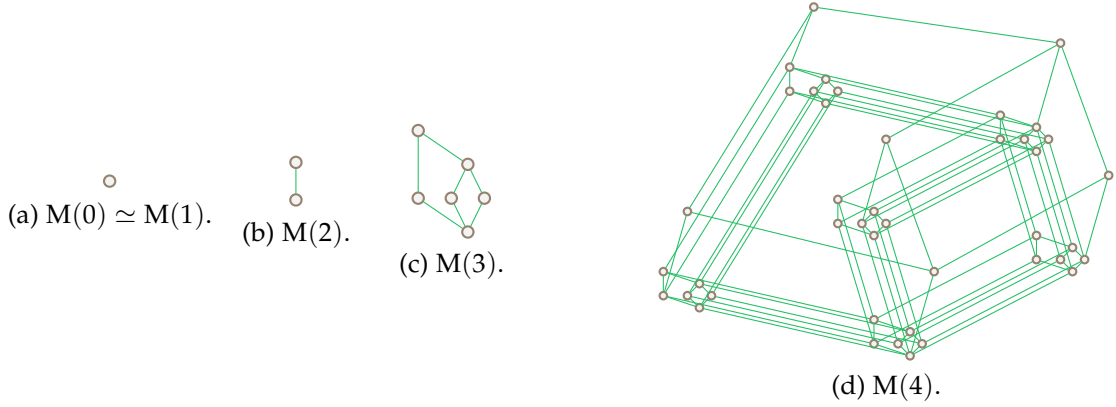

The \Def{Mockingbird lattice} of order $d \geq 0$ is the lattice $\MockingbirdLattice(d) :=
\Poset\Par{\RightComb_d}$ where $\RightComb_d$ is the combinator defined by
$\RightComb_0 := \M$ and, for any $d \geq 1$, by $\RightComb_d := \M \RightComb_{d - 1}$.
Figure~\ref{fig:hasse_diagram_mockingbird_lattices} shows the Hasse diagrams of the first
Mockingbird lattices.
These lattices are not graded, not self-dual, and not semidistributive.

\begin{Theorem} \label{thm:universality_mockingbird_lattices}
    For any $\ForestF \in \SetDuplicativeForests$, the poset
    $\SetDuplicativeForests(\ForestF)$ is isomorphic to a maximal interval of a Mockingbird
    lattice.
\end{Theorem}

Theorem~\ref{thm:universality_mockingbird_lattices} justifies the fact that the study of the
Mockingbird lattices is universal enough because these lattices contain as maximal interval
all duplicative forests lattices.

%%%%%%%%%%%%%%%%%%%%%%%%%%%%%%%%%%%%%%%%%%%%%%%%%%%%%%%%%%%%%%%%%%%%%%%%%%%%%%%%%%%%%%%%%%%%
%%%%%%%%%%%%%%%%%%%%%%%%%%%%%%%%%%%%%%%%%%%%%%%%%%%%%%%%%%%%%%%%%%%%%%%%%%%%%%%%%%%%%%%%%%%%
%%%%%%%%%%%%%%%%%%%%%%%%%%%%%%%%%%%%%%%%%%%%%%%%%%%%%%%%%%%%%%%%%%%%%%%%%%%%%%%%%%%%%%%%%%%%
\section{Enumerative properties} \label{sec:enumerative_properties}
Let $\K$ be any field of characteristic zero. For any set $X$, let $\K \Angle{X}$ be the
linear span of $X$. The dual space of $\K \Angle{X}$ is denoted by $\K \AAngle{X}$ and is by
definition the space of the maps $\SeriesF : X \to \K$, called \Def{$X$-series}. The
coefficient $\SeriesF(x)$ of any $x \in X$ is denoted by $\Angle{x, \SeriesF}$. The
\Def{support} of $\SeriesF$ is the set $\Support(\SeriesF) := \Bra{x \in X : \Angle{x,
\SeriesF} \ne 0}$. The \Def{characteristic series} of any subset $X'$ of $X$ is the series
$\CharacteristicSeries\Par{X'}$ having $X'$ as support and such that the coefficient of each
$x \in X'$ is $1$. For any $k \geq 0$, $\T^k \K \AAngle{X}$ is the $k$-th tensor power of
$\K \AAngle{X}$. Elements of this space are possibly infinite linear combinations of tensors
$x_1 \otimes \dots \otimes x_k$, where for any $i \in [k]$, $x_i \in X$. The \Def{tensor
algebra} of $\K \AAngle{X}$ is the space
\begin{math}
    \T^* \K \AAngle{X} := \bigoplus_{k \geq 0} \T^k \K \AAngle{X}.
\end{math}

A linear map $\phi : \T^{k_1} \K \AAngle{X} \to \T^{k_2} \K \AAngle{X}$, $k_1, k_2 \geq 0$,
is a \Def{$\Par{k_1, k_2}$-operation} on $\K \AAngle{X}$. The \Def{diagonal coproduct} is
the $(1, 2)$-operation $\Delta$ on $\K \AAngle{X}$ satisfying $\Delta(x) = x \otimes x$ for
any $x \in X$. When $X$ is endowed with an associative operation $\Product : X^2 \to X$, the
\Def{$\Product$-flattening map} is for any $k \geq 1$ the $(k, 1)$-operation
$\Flat^k_{\Product}$ on $\K \AAngle{X}$ satisfying
\begin{math}
    \Flat^k_{\Product}\Par{x_1 \otimes \dots \otimes x_k}
    = x_1 \Product \cdots \Product x_k
\end{math}
for any $x_1, \dots, x_k \in X$. When $X$ is endowed with an $n$-ary operation $\Product :
X^n \to X$, $n \geq 0$, the \Def{linearization} of $\Product$ is the $(n, 1)$-operation
$\bar{\Product}$ on $\K \AAngle{X}$ satisfying
\begin{math}
    \bar{\Product}\Par{x_1 \otimes \dots \otimes x_n} = \Product\Par{x_1, \dots, x_n}
\end{math}
for any $x_1, \dots, x_n \in X$. When $n = 1$, by a slight abuse of notation, for any $k
\geq 1$ and $x_1, \dots, x_k \in X$, we set
\begin{math}
    \bar{\Product}\Par{x_1 \otimes \dots \otimes x_k}
    := \bar{\Product}\Par{x_1} \otimes \dots \otimes \bar{\Product}\Par{x_k}.
\end{math}
To lighten the notation, when $\Product$ is a $(2, 1)$-operation on $\K \AAngle{X}$, we will
use $\Product$ as an infix operation by writing $\SeriesF_1 \Product \SeriesF_2$ for
$\Product\Par{\SeriesF_1 \otimes \SeriesF_2}$ for any $\SeriesF_1, \SeriesF_2 \in \K
\AAngle{X}$.

The space of the usual power series on the formal parameter $\VarZ$ is denoted by $\K
\AAngle{\VarZ}$. For any $F, F' \in \K \AAngle{\VarZ}$, $F\Han{\VarZ := F'}$ is the series
of $\K \AAngle{\VarZ}$ obtained by substituting $F'$ for $\VarZ$ in $F$. The \Def{Hadamard
product} is the binary operation $\HadamardProduct$ on $\K \AAngle{\VarZ}$ defined linearly
for any $n_1, n_2 \geq 0$ by
\begin{math}
    \VarZ^{n_1} \HadamardProduct \VarZ^{n_2}
    := \IndicatorFunction_{n_1 = n_2} \, \VarZ^{n_1}.
\end{math}
The \Def{max product} is the binary operation $\MaxProduct$ on $\K \AAngle{\VarZ}$ defined
linearly for any $n_1, n_2 \geq 0$ by $\VarZ^{n_1} \MaxProduct \VarZ^{n_2} := \VarZ^{\max
\Bra{n_1, n_2}}$. If $X$ is endowed with a map $\omega : X \to \N$, the
\Def{$\omega$-enumeration map} is the partial map $\EnumerationMap_\omega : \T^* \K
\AAngle{X} \to \K \AAngle{\VarZ}$ defined linearly for any $k \geq 1$ and $x_1, \dots, x_k
\in X$ by
\begin{math}
    \EnumerationMap_\omega\Par{x_1 \otimes \dots \otimes x_k}
    := \VarZ^{\omega\Par{x_1}} \MaxProduct \dots \MaxProduct \VarZ^{\omega\Par{x_k}}.
\end{math}
For any $\SeriesF \in \T^* \K \AAngle{X}$, the generating series
$\EnumerationMap_\omega(\SeriesF)$ is the \Def{$\omega$-enumeration} of $\SeriesF$. In the
sequel, we shall use the following strategy to enumerate a set $X$ w.r.t.\ such a map
$\omega$: we shall provide a description of $\CharacteristicSeries(X)$, then deduce a
description of $\EnumerationMap_\omega(\CharacteristicSeries(X))$, and finally deduce from
this a formula to compute the coefficients $\Angle{\VarZ^n,
\EnumerationMap_\omega(\CharacteristicSeries(X))}$, $n \geq 0$.

Recall now that $\Product$ is the binary operation on $\SetTerms(\GeneratingSet)$
satisfying, for any $\TreeT_1, \TreeT_2 \in \SetTerms(\GeneratingSet)$, $\TreeT_1 \Product
\TreeT_2 = \TreeT_1 \TreeT_2$. Proposition~\ref{prop:maximal_minimal} leads to the following
result.

\begin{Proposition} \label{prop:maximal_minimal_elements_enumeration}
    The characteristic series $\SeriesF_{\max}$ of the maximal combinators of $\Poset$
    satisfies
    \begin{equation} \label{equ:maximal_elements_enumeration}
        \SeriesF_{\max}
        = \M + \M \M + \SeriesF_{\max} \, \bar{\Product} \, \SeriesF_{\max}
        - \M \, \bar{\Product} \, \SeriesF_{\max}
    \end{equation}
    and the characteristic series $\SeriesF_{\min}$ of the minimal combinators of $\Poset$
    satisfies
    \begin{equation} \label{equ:minimal_elements_enumeration}
        \SeriesF_{\min}
        = \M + \M \M + \SeriesF_{\min} \, \bar{\Product} \, \SeriesF_{\min}
        - \bar{\Product}\Par{\Delta\Par{\SeriesF_{\min}}}.
    \end{equation}
\end{Proposition}

A consequence of Proposition~\ref{prop:maximal_minimal_elements_enumeration} is that the
$\Deg$-enumeration $F_{\max}$ of $\SeriesF_{\max}$ satisfies
\begin{math}
    F_{\max} = 1 + \VarZ + \VarZ F_{\max}^2 - \VarZ F_{\max}.
\end{math}
The first coefficients are $1$, $1$, $1$, $2$, $4$, $9$, $21$, $51$ and form
Sequence~\OEIS{A001006} (Motzkin numbers) of~\cite{Slo}. Another consequence of
Proposition~\ref{prop:maximal_minimal_elements_enumeration} is that the $\Deg$-enumeration
$F_{\min}$ of $\SeriesF_{\min}$ satisfies
\begin{math}
    F_{\min} = 1 + \VarZ + \VarZ F_{\min}^2 - \VarZ \, F_{\min}\Han{\VarZ := \VarZ^2}.
\end{math}
We deduce from this that the number of these terms of degree $d \geq 0$ is $\SequenceA(d)$
where $\SequenceA$ is the integer sequence satisfying $\SequenceA(0) = \SequenceA(1) = 1$
and, for any $d \geq 2$,
\begin{equation}
    \SequenceA(d) =
    \SequenceB(d - 1) - \IndicatorFunction_{d \mbox{ \scriptsize is odd}} \enspace
    \SequenceA\Par{(d - 1) / 2}
    \qquad \mbox{ where } \quad
    \SequenceB(d) := \sum_{i \in \HanL{d}} \SequenceA(i) \ \SequenceA(d - i).
\end{equation}
The first numbers are $1$, $1$, $2$, $4$, $12$, $34$, $108$, $344$ and form
Sequence~\OEIS{A343663} of~\cite{Slo}.

By Proposition~\ref{prop:first_graph_properties}, $\CLS$ has the properties described at the
very end of Section~\ref{sec:terms_rewrite_cls}. Therefore, the $\Height$-enumerations of
$\SeriesF_{\max}$ and $\SeriesF_{\min}$ are equal and is the generating series of the
$\Equiv$-equivalence classes of terms w.r.t.\ the height of their terms. By
Proposition~\ref{prop:maximal_minimal_elements_enumeration}, by denoting it by $F$, it
satisfies
\begin{math}
    F = 1 + \VarZ + \VarZ \Par{F \MaxProduct F} - \VarZ F.
\end{math}
Therefore, the number of these $\Equiv$-equivalence classes of terms of height $h \geq 0$ is
$\SequenceA(h)$ where $\SequenceA$ is the integer sequence satisfying $\SequenceA(0) =
\SequenceA(1) = 1$ and, for any $h \geq 2$,
\begin{equation}
    \SequenceA(h) =
        \SequenceA(h - 1)^2 - \SequenceA(h - 1)
        + 2 \SequenceA(h - 1) \sum_{i \in [h - 1]} \SequenceA(i - 1).
\end{equation}
The first numbers are $1$, $1$, $2$, $10$, $170$, $33490$, $1133870930$,
$1285739648704587610$ and form Sequence~\OEIS{A063573} of~\cite{Slo}.

Let us now consider series on duplicative forests in order to obtain enumerative results on
the Mockingbird lattices by using Proposition~\ref{prop:poset_isomorphism}. For any $k \geq
1$ and $u \in \Bra{\WhiteNode, \BlackNode}^k$, the \Def{merging product} is the $\Par{k +
|u|_{\BlackNode}, k}$-operation on $\K \AAngle{\SetDuplicativeForests}$ satisfying, for
any $\ForestF_1, \dots, \ForestF_{k + |u|_{\BlackNode}} \in \SetDuplicativeForests$,
\begin{math}
    \Merge_{\WhiteNode}\Par{\ForestF_1} = \WhiteNode \Par{\ForestF_1},
\end{math}
\begin{math}
        \Merge_{\WhiteNode u'}\Par{\ForestF_1 \otimes \dots \otimes \ForestF_k}
        = \Merge_{\WhiteNode}\Par{\ForestF_1}
            \otimes \Merge_{u'}\Par{\ForestF_2 \otimes \dots \otimes \ForestF_k},
\end{math}
\begin{math}
    \Merge_{\BlackNode}\Par{\ForestF_1 \otimes \ForestF_2}
    = \BlackNode \Par{\ForestF_1 \ForestF_2},
\end{math}
and
\begin{math}
    \Merge_{\BlackNode u'}\Par{\ForestF_1 \otimes \dots \otimes \ForestF_k}
    = \Merge_{\BlackNode}\Par{\ForestF_1 \otimes \ForestF_2}
        \otimes \Merge_{u'}\Par{\ForestF_3 \otimes \dots \otimes \ForestF_k},
\end{math}
where $u' \in \Bra{\WhiteNode, \BlackNode}^*$. For instance,
\begin{equation}
    \Merge_{\WhiteNode \, \BlackNode \, \BlackNode}
    \Par{
    \scalebox{.75}{
    \begin{tikzpicture}[Centering,xscale=0.17,yscale=0.3]
        \node[BlackNode](1)at(0.00,-1.00){};
        \node[BlackNode](0)at(0.00,0.00){};
        \draw[Edge](1)--(0);
    \end{tikzpicture}}
    \otimes
    \scalebox{.75}{
    \begin{tikzpicture}[Centering,xscale=0.17,yscale=0.2]
        \node[WhiteNode](0)at(0.00,-1.50){};
        \node[WhiteNode](2)at(2.00,-1.50){};
        \node[WhiteNode](1)at(1.00,0.00){};
        \draw[Edge](0)--(1);
        \draw[Edge](2)--(1);
    \end{tikzpicture}}
    \otimes
    \scalebox{.75}{
    \begin{tikzpicture}[Centering,xscale=0.17,yscale=0.24]
        \node[WhiteNode](0)at(0.00,-1.33){};
        \node[BlackNode](3)at(2.00,-2.67){};
        \node[WhiteNode](2)at(2.00,-1.33){};
        \draw[Edge](3)--(2);
    \end{tikzpicture}}
    \otimes
    \scalebox{.75}{
    \begin{tikzpicture}[Centering,xscale=0.14,yscale=0.2]
        \node[WhiteNode](0)at(0.00,-1.50){};
        \node[BlackNode](2)at(2.00,-1.50){};
    \end{tikzpicture}}
    \otimes
    \scalebox{.75}{
    \begin{tikzpicture}[Centering,xscale=0.14,yscale=0.2]
        \node[WhiteNode](0)at(0.00,-1.50){};
        \node[WhiteNode](2)at(2.00,-1.50){};
    \end{tikzpicture}}}
    =
    \scalebox{.75}{
    \begin{tikzpicture}[Centering,xscale=0.17,yscale=0.33]
        \node[BlackNode](2)at(0.00,-2.00){};
        \node[WhiteNode](0)at(0.00,0.00){};
        \node[BlackNode](1)at(0.00,-1.00){};
        \draw[Edge](1)--(0);
        \draw[Edge](2)--(1);
    \end{tikzpicture}}
    \otimes
    \scalebox{.75}{
    \begin{tikzpicture}[Centering,xscale=0.17,yscale=0.14]
        \node[WhiteNode](0)at(0.00,-4.67){};
        \node[WhiteNode](2)at(2.00,-4.67){};
        \node[WhiteNode](4)at(3.00,-2.33){};
        \node[BlackNode](6)at(5.00,-4.67){};
        \node[WhiteNode](1)at(1.00,-2.33){};
        \node[BlackNode](3)at(3.00,0.00){};
        \node[WhiteNode](5)at(5.00,-2.33){};
        \draw[Edge](0)--(1);
        \draw[Edge](1)--(3);
        \draw[Edge](2)--(1);
        \draw[Edge](4)--(3);
        \draw[Edge](5)--(3);
        \draw[Edge](6)--(5);
    \end{tikzpicture}}
    \otimes
    \scalebox{.75}{
    \begin{tikzpicture}[Centering,xscale=0.13,yscale=0.17]
        \node[WhiteNode](0)at(-1.00,-2.50){};
        \node[BlackNode](1)at(1.00,-2.50){};
        \node[WhiteNode](3)at(3.00,-2.50){};
        \node[WhiteNode](4)at(5.00,-2.50){};
        \node[BlackNode](2)at(2.00,0.00){};
        \draw[Edge](0)--(2);
        \draw[Edge](1)--(2);
        \draw[Edge](3)--(2);
        \draw[Edge](4)--(2);
    \end{tikzpicture}}.
\end{equation}

For any $d \geq 0$, the \Def{$d$-ladder} is the duplicative forest $\Ladder_d$ defined
recursively by $\Ladder_0 := \epsilon$ and, for any $d \geq 1$, by $\Ladder_d := \WhiteNode
\Par{\Ladder_{d - 1}}$. Let us denote by $\GreaterLadders$ the set $\bigcup_{d \geq 0}
\SetDuplicativeForests\Par{\Ladder_d}$. The \Def{series of ladders} is the unique
$\SetDuplicativeForests$-series $\SeriesLadders$ satisfying
\begin{math}
    \SeriesLadders = \epsilon + \bar{\WhiteNode}(\SeriesLadders).
\end{math}
Hence,
\begin{equation}
    \SeriesLadders = \sum_{d \geq 0} \Ladder_d
    =
    \epsilon
    +
    \scalebox{.75}{
    \begin{tikzpicture}[Centering,xscale=0.17,yscale=0.3]
        \node[WhiteNode](0)at(0.00,0.00){};
    \end{tikzpicture}}
    +
    \scalebox{.75}{
    \begin{tikzpicture}[Centering,xscale=0.17,yscale=0.3]
        \node[WhiteNode](1)at(0.00,-1.00){};
        \node[WhiteNode](0)at(0.00,0.00){};
    \draw[Edge](1)--(0);
    \end{tikzpicture}}
    +
    \scalebox{.75}{
    \begin{tikzpicture}[Centering,xscale=0.17,yscale=0.3]
        \node[WhiteNode](2)at(0.00,-2.00){};
        \node[WhiteNode](0)at(0.00,0.00){};
        \node[WhiteNode](1)at(0.00,-1.00){};
        \draw[Edge](1)--(0);
        \draw[Edge](2)--(1);
    \end{tikzpicture}}
    +
    \scalebox{.75}{
    \begin{tikzpicture}[Centering,xscale=0.17,yscale=0.3]
        \node[WhiteNode](3)at(0.00,-3.00){};
        \node[WhiteNode](0)at(0.00,0.00){};
        \node[WhiteNode](1)at(0.00,-1.00){};
        \node[WhiteNode](2)at(0.00,-2.00){};
        \draw[Edge](1)--(0);
        \draw[Edge](2)--(1);
        \draw[Edge](3)--(2);
    \end{tikzpicture}}
    + \cdots.
\end{equation}

Let $\SeriesGreater$ be the $(1, 1)$-operation on $\K \AAngle{\SetDuplicativeForests}$
satisfying, for any $\ForestF \in \SetDuplicativeForests$,
\begin{equation}
    \SeriesGreater(\ForestF)
    = \sum_{\ForestF' \in \SetDuplicativeForests(\ForestF)} \ForestF'.
\end{equation}
By definition, $\SeriesGreater(\ForestF)$ is the characteristic series of
$\SetDuplicativeForests(\ForestF)$. For instance,
\begin{equation}
    \SeriesGreater \Par{
    \scalebox{.75}{
    \begin{tikzpicture}[Centering,xscale=0.19,yscale=0.18]
        \node[WhiteNode](0)at(0.00,-4.00){};
        \node[BlackNode](3)at(2.00,-6.00){};
        \node[BlackNode](7)at(4.00,-6.00){};
        \node[BlackNode](1)at(1.00,-2.00){};
        \node[WhiteNode](2)at(2.00,-4.00){};
        \node[WhiteNode](5)at(4.00,-2.00){};
        \node[BlackNode](6)at(4.00,-4.00){};
        \draw[Edge](0)--(1);
        \draw[Edge](2)--(1);
        \draw[Edge](3)--(2);
        \draw[Edge](6)--(5);
        \draw[Edge](7)--(6);
    \end{tikzpicture}}}
    =
    \scalebox{.75}{
    \begin{tikzpicture}[Centering,xscale=0.2,yscale=0.18]
        \node[WhiteNode](0)at(0.00,-4.00){};
        \node[BlackNode](3)at(2.00,-6.00){};
        \node[BlackNode](7)at(4.00,-6.00){};
        \node[BlackNode](1)at(1.00,-2.00){};
        \node[WhiteNode](2)at(2.00,-4.00){};
        \node[WhiteNode](5)at(4.00,-2.00){};
        \node[BlackNode](6)at(4.00,-4.00){};
        \draw[Edge](0)--(1);
        \draw[Edge](2)--(1);
        \draw[Edge](3)--(2);
        \draw[Edge](6)--(5);
        \draw[Edge](7)--(6);
    \end{tikzpicture}}
    +
    \scalebox{.75}{
    \begin{tikzpicture}[Centering,xscale=0.2,yscale=0.18]
        \node[BlackNode](0)at(0.00,-4.00){};
        \node[BlackNode](3)at(2.00,-6.00){};
        \node[BlackNode](7)at(4.00,-6.00){};
        \node[BlackNode](1)at(1.00,-2.00){};
        \node[WhiteNode](2)at(2.00,-4.00){};
        \node[WhiteNode](5)at(4.00,-2.00){};
        \node[BlackNode](6)at(4.00,-4.00){};
        \draw[Edge](0)--(1);
        \draw[Edge](2)--(1);
        \draw[Edge](3)--(2);
        \draw[Edge](6)--(5);
        \draw[Edge](7)--(6);
    \end{tikzpicture}}
    +
    \scalebox{.75}{
    \begin{tikzpicture}[Centering,xscale=0.2,yscale=0.18]
        \node[WhiteNode](0)at(0.00,-4.00){};
        \node[BlackNode](3)at(2.75,-6.00){};
        \node[BlackNode](3')at(1.25,-6.00){};
        \node[BlackNode](7)at(5.00,-6.00){};
        \node[BlackNode](1)at(1.00,-2.00){};
        \node[BlackNode](2)at(2.00,-4.00){};
        \node[WhiteNode](5)at(5.00,-2.00){};
        \node[BlackNode](6)at(5.00,-4.00){};
        \draw[Edge](0)--(1);
        \draw[Edge](2)--(1);
        \draw[Edge](3)--(2);
        \draw[Edge](3')--(2);
        \draw[Edge](6)--(5);
        \draw[Edge](7)--(6);
    \end{tikzpicture}}
    +
    \scalebox{.75}{
    \begin{tikzpicture}[Centering,xscale=0.2,yscale=0.18]
        \node[WhiteNode](0)at(0.00,-4.00){};
        \node[BlackNode](3)at(2.00,-6.00){};
        \node[BlackNode](7)at(5.75,-6.00){};
        \node[BlackNode](7')at(4.25,-6.00){};
        \node[BlackNode](1)at(1.00,-2.00){};
        \node[WhiteNode](2)at(2.00,-4.00){};
        \node[BlackNode](5)at(5.00,-2.00){};
        \node[BlackNode](6)at(5.75,-4.00){};
        \node[BlackNode](6')at(4.25,-4.00){};
        \draw[Edge](0)--(1);
        \draw[Edge](2)--(1);
        \draw[Edge](3)--(2);
        \draw[Edge](6)--(5);
        \draw[Edge](6')--(5);
        \draw[Edge](7)--(6);
        \draw[Edge](7')--(6');
    \end{tikzpicture}}
    +
    \scalebox{.75}{
    \begin{tikzpicture}[Centering,xscale=0.2,yscale=0.18]
        \node[BlackNode](0)at(0.00,-4.00){};
        \node[BlackNode](3)at(2.75,-6.00){};
        \node[BlackNode](3')at(1.25,-6.00){};
        \node[BlackNode](7)at(5.00,-6.00){};
        \node[BlackNode](1)at(1.00,-2.00){};
        \node[BlackNode](2)at(2.00,-4.00){};
        \node[WhiteNode](5)at(5.00,-2.00){};
        \node[BlackNode](6)at(5.00,-4.00){};
        \draw[Edge](0)--(1);
        \draw[Edge](2)--(1);
        \draw[Edge](3)--(2);
        \draw[Edge](3')--(2);
        \draw[Edge](6)--(5);
        \draw[Edge](7)--(6);
    \end{tikzpicture}}
    +
    \scalebox{.75}{
    \begin{tikzpicture}[Centering,xscale=0.2,yscale=0.18]
        \node[BlackNode](0)at(0.00,-4.00){};
        \node[BlackNode](3)at(2.00,-6.00){};
        \node[BlackNode](7)at(5.75,-6.00){};
        \node[BlackNode](7')at(4.25,-6.00){};
        \node[BlackNode](1)at(1.00,-2.00){};
        \node[WhiteNode](2)at(2.00,-4.00){};
        \node[BlackNode](5)at(5.00,-2.00){};
        \node[BlackNode](6)at(5.75,-4.00){};
        \node[BlackNode](6')at(4.25,-4.00){};
        \draw[Edge](0)--(1);
        \draw[Edge](2)--(1);
        \draw[Edge](3)--(2);
        \draw[Edge](6)--(5);
        \draw[Edge](6')--(5);
        \draw[Edge](7)--(6);
        \draw[Edge](7')--(6');
    \end{tikzpicture}}
    +
    \scalebox{.75}{
    \begin{tikzpicture}[Centering,xscale=0.2,yscale=0.18]
        \node[WhiteNode](0)at(0.00,-4.00){};
        \node[BlackNode](3)at(2.75,-6.00){};
        \node[BlackNode](3')at(1.25,-6.00){};
        \node[BlackNode](7)at(5.75,-6.00){};
        \node[BlackNode](7')at(4.25,-6.00){};
        \node[BlackNode](1)at(1.00,-2.00){};
        \node[BlackNode](2)at(2.00,-4.00){};
        \node[BlackNode](5)at(5.00,-2.00){};
        \node[BlackNode](6)at(5.75,-4.00){};
        \node[BlackNode](6')at(4.25,-4.00){};
        \draw[Edge](0)--(1);
        \draw[Edge](2)--(1);
        \draw[Edge](3)--(2);
        \draw[Edge](3')--(2);
        \draw[Edge](6)--(5);
        \draw[Edge](6')--(5);
        \draw[Edge](7)--(6);
        \draw[Edge](7')--(6');
    \end{tikzpicture}}
    +
    \scalebox{.75}{
    \begin{tikzpicture}[Centering,xscale=0.2,yscale=0.18]
        \node[BlackNode](0)at(0.00,-4.00){};
        \node[BlackNode](3)at(2.75,-6.00){};
        \node[BlackNode](3')at(1.25,-6.00){};
        \node[BlackNode](7)at(5.75,-6.00){};
        \node[BlackNode](7')at(4.25,-6.00){};
        \node[BlackNode](1)at(1.00,-2.00){};
        \node[BlackNode](2)at(2.00,-4.00){};
        \node[BlackNode](5)at(5.00,-2.00){};
        \node[BlackNode](6)at(5.75,-4.00){};
        \node[BlackNode](6')at(4.25,-4.00){};
        \draw[Edge](0)--(1);
        \draw[Edge](2)--(1);
        \draw[Edge](3)--(2);
        \draw[Edge](3')--(2);
        \draw[Edge](6)--(5);
        \draw[Edge](6')--(5);
        \draw[Edge](7)--(6);
        \draw[Edge](7')--(6');
    \end{tikzpicture}}.
\end{equation}
Observe that $\SeriesGreater(\SeriesLadders)$ is the characteristic series of
$\GreaterLadders$ and that $\EnumerationMap_\Height(\SeriesGreater(\SeriesLadders))$ is the
generating series of the cardinalities of the lattices
$\SetDuplicativeForests\Par{\Ladder_d}$, enumerated w.r.t.\ $d \geq 0$.

\begin{Theorem} \label{thm:series_elements}
    The series $\SeriesGreater(\SeriesLadders)$ satisfies
    \begin{equation}
        \SeriesGreater(\SeriesLadders)
        = \epsilon + \bar{\WhiteNode}(\SeriesGreater(\SeriesLadders))
        + \bar{\BlackNode}\Par{\SeriesGreater\Par{\Flat^2_{\ConcatenateForests}
            \Par{\Delta\Par{\SeriesLadders}}}}.
    \end{equation}
\end{Theorem}

We deduce from Theorem~\ref{thm:series_elements} that the $\Height$-enumeration $F$ of
$\SeriesGreater(\SeriesLadders)$ satisfies
\begin{math}
    F = 1 + \VarZ F + \VarZ (F \HadamardProduct F)
\end{math}
so that for any $d \geq 1$, the number of elements in $\MockingbirdLattice(d)$ is
$\SequenceA(d - 1)$ where $\SequenceA$ is the integer sequence satisfying $\SequenceA(0) =
1$ and, for any $d \geq 1$,
\begin{equation} \label{equ:recurrence_number_elements}
    \SequenceA(d) = \SequenceA(d - 1) + \SequenceA(d - 1)^2.
\end{equation}
The sequence of the cardinalities of $\MockingbirdLattice(d)$, $d \geq 0$, starts by $1$,
$1$, $2$, $6$, $42$, $1806$, $3263442$, $10650056950806$ and forms Sequence~\OEIS{A007018}
of~\cite{Slo}.

Let $\SeriesCoverings$ be the $(1, 1)$-operation on $\K \AAngle{\SetDuplicativeForests}$
satisfying, for any $\ForestF \in \SetDuplicativeForests$,
\begin{equation}
    \SeriesCoverings(\ForestF)
    = \sum_{\substack{
        \ForestF' \in \SetDuplicativeForests \\
        \ForestF \RewDuplicative \ForestF'
    }}
    \ForestF'.
\end{equation}
Let also $\SeriesNbInputs$ be the $(1, 1)$-operation on $\K \AAngle{\SetDuplicativeForests}$
satisfying $\SeriesNbInputs(\ForestF) = \SeriesCoverings(\SeriesGreater(\ForestF))$ for any
$\ForestF \in \SetDuplicativeForests$. By a straightforward computation, we obtain
\begin{equation}
    \SeriesNbInputs(\ForestF) =
    \sum_{\ForestF' \in \SetDuplicativeForests(\ForestF)}
    \# \Bra{\ForestF'' \in \SetDuplicativeForests(\ForestF) : \ForestF'' \RewDuplicative
    \ForestF'} \, \ForestF',
\end{equation}
so that the coefficient of each $\ForestF' \in \SetDuplicativeForests(\ForestF)$ in
$\SeriesNbInputs(\ForestF)$ is the number of duplicative forests admitting $\ForestF'$ as
covering in $\SetDuplicativeForests(\ForestF)$. For instance (see at the same time
Figure~\ref{fig:example_poset_duplicative_forests}),
\begin{multline}
    \SeriesNbInputs\Par{
    \scalebox{.75}{
    \begin{tikzpicture}[Centering,xscale=0.18,yscale=0.26]
        \node[WhiteNode](1)at(0.00,-2.67){};
        \node[WhiteNode](3)at(2.00,-1.33){};
        \node[WhiteNode](0)at(0.00,-1.33){};
        \draw[Edge](1)--(0);
    \end{tikzpicture}}}
    =
    \scalebox{.75}{
    \begin{tikzpicture}[Centering,xscale=0.18,yscale=0.26]
        \node[WhiteNode](1)at(0.00,-2.67){};
        \node[BlackNode](3)at(2.00,-1.33){};
        \node[WhiteNode](0)at(0.00,-1.33){};
        \draw[Edge](1)--(0);
    \end{tikzpicture}}
    +
    \scalebox{.75}{
    \begin{tikzpicture}[Centering,xscale=0.18,yscale=0.26]
        \node[BlackNode](1)at(0.00,-2.67){};
        \node[WhiteNode](3)at(2.00,-1.33){};
        \node[WhiteNode](0)at(0.00,-1.33){};
        \draw[Edge](1)--(0);
    \end{tikzpicture}}
    +
    2
    \scalebox{.75}{
    \begin{tikzpicture}[Centering,xscale=0.18,yscale=0.26]
        \node[BlackNode](1)at(0.00,-2.67){};
        \node[BlackNode](3)at(2.00,-1.33){};
        \node[WhiteNode](0)at(0.00,-1.33){};
        \draw[Edge](1)--(0);
    \end{tikzpicture}}
    +
    \scalebox{.75}{
    \begin{tikzpicture}[Centering,xscale=0.18,yscale=0.21]
        \node[WhiteNode](0)at(0.00,-3.33){};
        \node[WhiteNode](2)at(2.00,-3.33){};
        \node[WhiteNode](4)at(4.00,-1.67){};
        \node[BlackNode](1)at(1.00,-1.67){};
        \draw[Edge](0)--(1);
        \draw[Edge](2)--(1);
    \end{tikzpicture}}
    +
    2
    \scalebox{.75}{
    \begin{tikzpicture}[Centering,xscale=0.18,yscale=0.21]
        \node[WhiteNode](0)at(0.00,-3.33){};
        \node[WhiteNode](2)at(2.00,-3.33){};
        \node[BlackNode](4)at(4.00,-1.67){};
        \node[BlackNode](1)at(1.00,-1.67){};
        \draw[Edge](0)--(1);
        \draw[Edge](2)--(1);
    \end{tikzpicture}}
    +
    \scalebox{.75}{
    \begin{tikzpicture}[Centering,xscale=0.18,yscale=0.21]
        \node[WhiteNode](0)at(0.00,-3.33){};
        \node[BlackNode](2)at(2.00,-3.33){};
        \node[WhiteNode](4)at(4.00,-1.67){};
        \node[BlackNode](1)at(1.00,-1.67){};
        \draw[Edge](0)--(1);
        \draw[Edge](2)--(1);
    \end{tikzpicture}}
    +
    \scalebox{.75}{
    \begin{tikzpicture}[Centering,xscale=0.18,yscale=0.21]
        \node[BlackNode](0)at(0.00,-3.33){};
        \node[WhiteNode](2)at(2.00,-3.33){};
        \node[WhiteNode](4)at(4.00,-1.67){};
        \node[BlackNode](1)at(1.00,-1.67){};
        \draw[Edge](0)--(1);
        \draw[Edge](2)--(1);
    \end{tikzpicture}}
    +
    2
    \scalebox{.75}{
    \begin{tikzpicture}[Centering,xscale=0.18,yscale=0.21]
        \node[WhiteNode](0)at(0.00,-3.33){};
        \node[BlackNode](2)at(2.00,-3.33){};
        \node[BlackNode](4)at(4.00,-1.67){};
        \node[BlackNode](1)at(1.00,-1.67){};
        \draw[Edge](0)--(1);
        \draw[Edge](2)--(1);
    \end{tikzpicture}}
    +
    2
    \scalebox{.75}{
    \begin{tikzpicture}[Centering,xscale=0.18,yscale=0.21]
        \node[BlackNode](0)at(0.00,-3.33){};
        \node[WhiteNode](2)at(2.00,-3.33){};
        \node[BlackNode](4)at(4.00,-1.67){};
        \node[BlackNode](1)at(1.00,-1.67){};
        \draw[Edge](0)--(1);
        \draw[Edge](2)--(1);
    \end{tikzpicture}}
    +
    3
    \scalebox{.75}{
    \begin{tikzpicture}[Centering,xscale=0.18,yscale=0.21]
        \node[BlackNode](0)at(0.00,-3.33){};
        \node[BlackNode](2)at(2.00,-3.33){};
        \node[WhiteNode](4)at(4.00,-1.67){};
        \node[BlackNode](1)at(1.00,-1.67){};
        \draw[Edge](0)--(1);
        \draw[Edge](2)--(1);
    \end{tikzpicture}}
    +
    4
    \scalebox{.75}{
    \begin{tikzpicture}[Centering,xscale=0.18,yscale=0.21]
        \node[BlackNode](0)at(0.00,-3.33){};
        \node[BlackNode](2)at(2.00,-3.33){};
        \node[BlackNode](4)at(4.00,-1.67){};
        \node[BlackNode](1)at(1.00,-1.67){};
        \draw[Edge](0)--(1);
        \draw[Edge](2)--(1);
    \end{tikzpicture}}.
\end{multline}
Observe that
\begin{math}
    \Support(\SeriesNbInputs(\SeriesLadders))
    = \GreaterLadders \setminus \Bra{\Ladder_d : d \geq 0}
\end{math}
and that $\EnumerationMap_{\Height}(\SeriesNbInputs(\SeriesLadders))$ is the generating
series of the number of edges of the Hasse diagrams of the lattices
$\SetDuplicativeForests\Par{\Ladder_d}$, enumerated w.r.t.\ $d \geq 0$.

\begin{Theorem} \label{thm:series_number_coverings}
    The series $\SeriesNbInputs(\SeriesLadders)$ satisfies
    \begin{equation}
        \SeriesNbInputs(\SeriesLadders)
        = \bar{\WhiteNode}(\SeriesNbInputs(\SeriesLadders))
        + \bar{\BlackNode}\Par{\SeriesNbInputs\Par{\Flat^2_{\ConcatenateForests}\Par{
            \Delta\Par{\SeriesLadders}}}}
        + \bar{\BlackNode}\Par{\Flat^2_{\ConcatenateForests}\Par{\Delta\Par{
            \SeriesGreater\Par{\SeriesLadders}}}}.
    \end{equation}
\end{Theorem}

We deduce from Theorem~\ref{thm:series_number_coverings} that the $\Height$-enumeration $F$
of $\SeriesNbInputs(\SeriesLadders)$ satisfies
\begin{math}
    F = \VarZ F + \VarZ G + 2 \VarZ (F \HadamardProduct G)
\end{math}
where $G$ is the $\Height$-enumeration of $\SeriesGreater(\SeriesLadders)$. Therefore, for
any $d \geq 1$, the number of edges in the Hasse diagram of $\MockingbirdLattice(d)$ is
$\SequenceA(d - 1)$ where $\SequenceA$ is the integer sequence satisfying $\SequenceA(0) =
0$ and, for any $d \geq 1$,
\begin{equation}
    \SequenceA(d)
    =  \SequenceA(d - 1) + \SequenceB(d - 1) + 2 \SequenceA(d - 1) \SequenceB(d - 1),
\end{equation}
where $\SequenceB$ is the integer sequence such that for any $d \geq 0$, $\SequenceB(d)$ is
the number of elements of $\SetDuplicativeForests\Par{\Ladder_d}$, satisfying
therefore~\eqref{equ:recurrence_number_elements}. The sequence of the number of edges of the
Hasse diagram of $\MockingbirdLattice(d)$, $d \geq 0$, starts by $0$, $0$, $1$, $7$, $97$,
$8287$, $29942737$, $195432804247687$. This sequence does not appear in~\cite{Slo} for the
time being.

Now, let $\SeriesNbSmaller$ be the $(1, 1)$-operation on $\K
\AAngle{\SetDuplicativeForests}$ satisfying $\SeriesNbSmaller(\ForestF) =
\SeriesGreater(\SeriesGreater(\ForestF))$ for any $\ForestF \in \SetDuplicativeForests$. By
a straightforward computation, we obtain
\begin{equation}
    \SeriesNbSmaller(\ForestF) =
    \sum_{\ForestF' \in \SetDuplicativeForests(\ForestF)}
    \# \Han{\ForestF, \ForestF'} \ForestF',
\end{equation}
so that the coefficient of each $\ForestF' \in \SetDuplicativeForests(\ForestF)$ in
$\SeriesNbSmaller(\ForestF)$ is the number of duplicative forests smaller than or equal to
$\ForestF'$ in $\SetDuplicativeForests(\ForestF)$. For instance (see at the same time
Figure~\ref{fig:example_poset_duplicative_forests}),
\begin{multline}
    \SeriesNbSmaller\Par{
    \scalebox{.75}{
    \begin{tikzpicture}[Centering,xscale=0.18,yscale=0.26]
        \node[WhiteNode](1)at(0.00,-2.67){};
        \node[WhiteNode](3)at(2.00,-1.33){};
        \node[WhiteNode](0)at(0.00,-1.33){};
        \draw[Edge](1)--(0);
    \end{tikzpicture}}}
    =
    \scalebox{.75}{
    \begin{tikzpicture}[Centering,xscale=0.18,yscale=0.26]
        \node[WhiteNode](1)at(0.00,-2.67){};
        \node[WhiteNode](3)at(2.00,-1.33){};
        \node[WhiteNode](0)at(0.00,-1.33){};
        \draw[Edge](1)--(0);
    \end{tikzpicture}}
    +
    2
    \scalebox{.75}{
    \begin{tikzpicture}[Centering,xscale=0.18,yscale=0.26]
        \node[WhiteNode](1)at(0.00,-2.67){};
        \node[BlackNode](3)at(2.00,-1.33){};
        \node[WhiteNode](0)at(0.00,-1.33){};
        \draw[Edge](1)--(0);
    \end{tikzpicture}}
    +
    2
    \scalebox{.75}{
    \begin{tikzpicture}[Centering,xscale=0.18,yscale=0.26]
        \node[BlackNode](1)at(0.00,-2.67){};
        \node[WhiteNode](3)at(2.00,-1.33){};
        \node[WhiteNode](0)at(0.00,-1.33){};
        \draw[Edge](1)--(0);
    \end{tikzpicture}}
    +
    4
    \scalebox{.75}{
    \begin{tikzpicture}[Centering,xscale=0.18,yscale=0.26]
        \node[BlackNode](1)at(0.00,-2.67){};
        \node[BlackNode](3)at(2.00,-1.33){};
        \node[WhiteNode](0)at(0.00,-1.33){};
        \draw[Edge](1)--(0);
    \end{tikzpicture}}
    +
    2
    \scalebox{.75}{
    \begin{tikzpicture}[Centering,xscale=0.18,yscale=0.21]
        \node[WhiteNode](0)at(0.00,-3.33){};
        \node[WhiteNode](2)at(2.00,-3.33){};
        \node[WhiteNode](4)at(4.00,-1.67){};
        \node[BlackNode](1)at(1.00,-1.67){};
        \draw[Edge](0)--(1);
        \draw[Edge](2)--(1);
    \end{tikzpicture}}
    +
    4
    \scalebox{.75}{
    \begin{tikzpicture}[Centering,xscale=0.18,yscale=0.21]
        \node[WhiteNode](0)at(0.00,-3.33){};
        \node[WhiteNode](2)at(2.00,-3.33){};
        \node[BlackNode](4)at(4.00,-1.67){};
        \node[BlackNode](1)at(1.00,-1.67){};
        \draw[Edge](0)--(1);
        \draw[Edge](2)--(1);
    \end{tikzpicture}}
    +
    3
    \scalebox{.75}{
    \begin{tikzpicture}[Centering,xscale=0.18,yscale=0.21]
        \node[WhiteNode](0)at(0.00,-3.33){};
        \node[BlackNode](2)at(2.00,-3.33){};
        \node[WhiteNode](4)at(4.00,-1.67){};
        \node[BlackNode](1)at(1.00,-1.67){};
        \draw[Edge](0)--(1);
        \draw[Edge](2)--(1);
    \end{tikzpicture}}
    +
    3
    \scalebox{.75}{
    \begin{tikzpicture}[Centering,xscale=0.18,yscale=0.21]
        \node[BlackNode](0)at(0.00,-3.33){};
        \node[WhiteNode](2)at(2.00,-3.33){};
        \node[WhiteNode](4)at(4.00,-1.67){};
        \node[BlackNode](1)at(1.00,-1.67){};
        \draw[Edge](0)--(1);
        \draw[Edge](2)--(1);
    \end{tikzpicture}}
    +
    6
    \scalebox{.75}{
    \begin{tikzpicture}[Centering,xscale=0.18,yscale=0.21]
        \node[WhiteNode](0)at(0.00,-3.33){};
        \node[BlackNode](2)at(2.00,-3.33){};
        \node[BlackNode](4)at(4.00,-1.67){};
        \node[BlackNode](1)at(1.00,-1.67){};
        \draw[Edge](0)--(1);
        \draw[Edge](2)--(1);
    \end{tikzpicture}}
    +
    6
    \scalebox{.75}{
    \begin{tikzpicture}[Centering,xscale=0.18,yscale=0.21]
        \node[BlackNode](0)at(0.00,-3.33){};
        \node[WhiteNode](2)at(2.00,-3.33){};
        \node[BlackNode](4)at(4.00,-1.67){};
        \node[BlackNode](1)at(1.00,-1.67){};
        \draw[Edge](0)--(1);
        \draw[Edge](2)--(1);
    \end{tikzpicture}}
    \\
    +
    6
    \scalebox{.75}{
    \begin{tikzpicture}[Centering,xscale=0.18,yscale=0.21]
        \node[BlackNode](0)at(0.00,-3.33){};
        \node[BlackNode](2)at(2.00,-3.33){};
        \node[WhiteNode](4)at(4.00,-1.67){};
        \node[BlackNode](1)at(1.00,-1.67){};
        \draw[Edge](0)--(1);
        \draw[Edge](2)--(1);
    \end{tikzpicture}}
    +
    12
    \scalebox{.75}{
    \begin{tikzpicture}[Centering,xscale=0.18,yscale=0.21]
        \node[BlackNode](0)at(0.00,-3.33){};
        \node[BlackNode](2)at(2.00,-3.33){};
        \node[BlackNode](4)at(4.00,-1.67){};
        \node[BlackNode](1)at(1.00,-1.67){};
        \draw[Edge](0)--(1);
        \draw[Edge](2)--(1);
    \end{tikzpicture}}.
\end{multline}

Contrary to what we have undertaken previously to express $\SeriesGreater(\SeriesLadders)$
and $\SeriesNbInputs(\SeriesLadders)$, we fail to directly express
$\SeriesNbSmaller(\SeriesLadders)$. The trick here consists in considering first a slightly
different series depending on a parameter $k \geq 1$ which can be seen as a catalytic
parameter. For any $k \geq 1$, let $\SeriesMeetDecomposition_k$ be the $(1, k)$-operation on
$\K \AAngle{\SetDuplicativeForests}$ satisfying, for any $\ForestF \in
\SetDuplicativeForests$,
\begin{equation}
    \SeriesMeetDecomposition_k(\ForestF) =
    \sum_{\substack{
        \ForestG_1, \dots, \ForestG_k \in \SetDuplicativeForests(\ForestF) \\
        \ForestG_1 \Meet \dots \Meet \ForestG_k = \ForestF
    }}
    \ForestG_1 \otimes \dots \otimes \ForestG_k.
\end{equation}
We call $\SeriesMeetDecomposition_k(\ForestF)$ the \Def{meet $k$-decomposition} of
$\ForestF$. Observe that $\SeriesMeetDecomposition_1$ is the identity map.

Observe that $\Support(\SeriesNbSmaller(\SeriesLadders)) = \GreaterLadders$ and that
\begin{math}
    \EnumerationMap_{\Height}\Par{\SeriesMeetDecomposition_1\Par{
        \SeriesNbSmaller\Par{\SeriesLadders}}}
\end{math}
is the generating series of the number of intervals of the lattices
$\SetDuplicativeForests\Par{\Ladder_d}$, enumerated w.r.t.\ $d \geq 0$.

\begin{Theorem} \label{thm:series_nb_smaller}
    The series $\SeriesNbSmaller(\SeriesLadders)$ satisfies
    $\SeriesNbSmaller(\SeriesLadders) =
    \SeriesMeetDecomposition_1(\SeriesNbSmaller(\SeriesLadders))$ where, for any $k \geq 1$,
    the series $\SeriesMeetDecomposition_k(\SeriesNbSmaller(\SeriesLadders))$ satisfies
    \begin{equation}
        \SeriesMeetDecomposition_k(\SeriesNbSmaller(\SeriesLadders))
        = \epsilon^{\otimes k}
        + \sum_{u \in \Bra{\WhiteNode, \BlackNode}^k}
            \Merge_u\Par{\SeriesMeetDecomposition_{k + |u|_{\BlackNode}}
                (\SeriesNbSmaller(\SeriesLadders))}
        + \bar{\BlackNode} \Par{
            \SeriesMeetDecomposition_k\Par{\SeriesNbSmaller\Par{
                    \Flat^k_{\ConcatenateForests}
            \Par{\Delta\Par{\SeriesLadders}}}}}.
    \end{equation}
\end{Theorem}

We deduce from Theorem~\ref{thm:series_nb_smaller} that the $\Height$-enumeration $F$ of
$\SeriesNbSmaller(\SeriesLadders)$ satisfies $F = F_1$ where, for any $k \geq 1$, $F_k$ is
the $\Height$-enumeration of $\SeriesMeetDecomposition_k(\SeriesNbSmaller(\SeriesLadders))$
which satisfies
\begin{math}
    F_k = 1 + \VarZ \Par{F_k \HadamardProduct F_k}
        + \VarZ \sum_{i \in \HanL{k}} \binom{k}{i} F_{k + i}.
\end{math}
Therefore, for any $d \geq 1$, the number of intervals in $\MockingbirdLattice(d)$ is
$\SequenceA_1(d - 1)$ where for any $k \geq 1$, $\SequenceA_k$ is the integer sequence
satisfying $\SequenceA_k(0) = 1$ and, for any $d \geq 1$,
\begin{equation}
    \SequenceA_k(d) =
    \SequenceA_k(d - 1)^2
    + \sum_{i \in \HanL{k}} \binom{k}{i} \SequenceA_{k + i}(d - 1).
\end{equation}
The sequence of the number of intervals of $\MockingbirdLattice(d)$, $d \geq 0$, starts by
$1$, $1$, $3$, $17$, $371$, $144513$, $20932611523$, $438176621806663544657$. This sequence
does not appear in~\cite{Slo} for the time being.

%%%%%%%%%%%%%%%%%%%%%%%%%%%%%%%%%%%%%%%%%%%%%%%%%%%%%%%%%%%%%%%%%%%%%%%%%%%%%%%%%%%%%%%%%%%%
%%%%%%%%%%%%%%%%%%%%%%%%%%%%%%%%%%%%%%%%%%%%%%%%%%%%%%%%%%%%%%%%%%%%%%%%%%%%%%%%%%%%%%%%%%%%
%%%%%%%%%%%%%%%%%%%%%%%%%%%%%%%%%%%%%%%%%%%%%%%%%%%%%%%%%%%%%%%%%%%%%%%%%%%%%%%%%%%%%%%%%%%%
\section*{Open questions and future work}
We have studied a CLS having many rich combinatorial properties despite its simplicity. This
can be considered as the prototypical example for this kind of investigation. We expect to
discover similar properties for more complex CLS. Additionally, here are three open
questions raised by this work.
\begin{enumerate*}[label={\bf (\arabic*)}]
    \item The description of minimal and maximal elements of $\Poset$ uses a notion of
    pattern avoidance in terms. This is a general fact: when a CLS $(\GeneratingSet, \Rew)$
    has the poset property, its minimal (resp.\ maximal) elements are the terms avoiding
    terms deduced from the ones appearing as right-hand (resp.\ left-hand) members of
    $\Rew$. Such an enumerative problem has been considered in~\cite{Gir20} for the
    particular case of terms without repeating variables. We ask here for the general
    enumeration of terms avoiding a set of terms wherein multiple occurrences of a same
    variable are allowed.
    \item We have shown that the Mockingbird CLS has the poset property, is rooted, and has
    the lattice property by employing some specific reasoning from the definition of the
    basic combinator $\M$. A question here concerns the existence of a general criterion to
    decide if a CLS has the poset (resp.\ lattice) property and if it is rooted.
    \item Finally, we have seen from Proposition~\ref{prop:finite_equivalence_classes} that
    being hierarchical is a sufficient condition for a CLS $\CLS$ to be locally finite. The
    question in this context consists in strengthening this result in order to obtain a
    necessary and sufficient condition for this last property.
\end{enumerate*}

%%%%%%%%%%%%%%%%%%%%%%%%%%%%%%%%%%%%%%%%%%%%%%%%%%%%%%%%%%%%%%%%%%%%%%%%%%%%%%%%%%%%%%%%%%%%
%%%%%%%%%%%%%%%%%%%%%%%%%%%%%%%%%%%%%%%%%%%%%%%%%%%%%%%%%%%%%%%%%%%%%%%%%%%%%%%%%%%%%%%%%%%%
%%%%%%%%%%%%%%%%%%%%%%%%%%%%%%%%%%%%%%%%%%%%%%%%%%%%%%%%%%%%%%%%%%%%%%%%%%%%%%%%%%%%%%%%%%%%
\bibliographystyle{plain}
\bibliography{Bibliography}

\begin{thebibliography}{10}

\bibitem{BEJW17}
H.~Barendregt, J.~Endrullis, J.~W. Klop, and J.~Waldmann.
\newblock {Dance of the starlings}.
\newblock In {\em {Raymond Smullyan on self reference}}, volume~14 of {\em
  Outst. Contrib. Log.}, pages 67--111. Springer, Cham, 2017.

\bibitem{BKVT03}
M.~Bezem, J.~W. Klop, R.~de~Vrijer, and Terese.
\newblock {\em {Term Rewriting Systems}}.
\newblock Cambridge Tracts in Theoretical Computer Science. Cambridge
  University Press, 2003.

\bibitem{Cur30}
H.~B. Curry.
\newblock {Grundlagen der Kombinatorischen Logik}.
\newblock {\em Am. J. Math.}, 52(4):789--834, 1930.

\bibitem{Gir20}
S.~Giraudo.
\newblock {Tree series and pattern avoidance in syntax trees}.
\newblock {\em J. Comb. Theory A}, 176, 2020.

\bibitem{HS08}
J.~R. Hindley and J.~P. Seldin.
\newblock {\em {Lambda-calculus and combinators, an introduction}}.
\newblock Cambridge University Press, Cambridge, 2008.

\bibitem{PS01}
D.~Probst and T.~Studer.
\newblock {How to normalize the Jay}.
\newblock {\em Theoret. Comput. Sci.}, 254(1-2):677--681, 2001.

\bibitem{Sch24}
M.~Schönfinkel.
\newblock {Über die Bausteine der mathematischen Logik}.
\newblock {\em Mathematische Annalen}, 92:305--316, 1924.

\bibitem{Slo}
N.~J.~A. Sloane.
\newblock {The On-Line Encyclopedia of Integer Sequences}.
\newblock \url{https://oeis.org/}.

\bibitem{Smu85}
R.~Smullyan.
\newblock {\em {To Mock a Mockingbird}}.
\newblock Alfred A. Knopf, Inc., 1985.

\bibitem{SWB93}
M.~Sprenger and M.~Wymann-Böni.
\newblock {How to decide the lark}.
\newblock {\em Theoret. Comput. Sci.}, 110(2):419--432, 1993.

\bibitem{Sta89}
R.~Statman.
\newblock {The word problem for Smullyan's lark combinator is decidable}.
\newblock {\em J. Symbolic Comput.}, 7(2):103--112, 1989.

\bibitem{Sta00}
R.~Statman.
\newblock {On the Word Problem for Combinators}.
\newblock {\em Lect. Notes Comp. Sci., Rewriting Techniques and Applications,
  11th International Conference}, 1833:203--213, 2000.

\bibitem{Sta17}
R.~Statman.
\newblock {Some tweets about Mockingbirds}.
\newblock In {\em {Raymond Smullyan on self reference}}, volume~14 of {\em
  Outst. Contrib. Log.}, pages 113--122. Springer, Cham, 2017.

\bibitem{Tam62}
D.~Tamari.
\newblock {The algebra of bracketings and their enumeration}.
\newblock {\em Nieuw Arch. Wisk.}, 10(3):131--146, 1962.

\bibitem{Wal00}
J.~Waldmann.
\newblock {The combinator {${\bf S}$}}.
\newblock {\em Inform. Comput.}, 159(1-2):2--21, 2000.

\end{thebibliography}

\end{document}